\documentclass[11pt,oneside]{amsart}
\usepackage{bm}
\usepackage{amssymb,amsfonts}
\usepackage{amsmath,latexsym}
\usepackage{xcolor}
\usepackage{mathrsfs}
\usepackage{enumerate}

\usepackage[top=1in, bottom=1in, left=1.25in, right=1.25in]{geometry}

\usepackage[unicode,breaklinks=true,colorlinks=true]{hyperref}

\newtheorem{thm}{Theorem}[section]

\newtheorem{lem}[thm]{Lemma}

\theoremstyle{definition}
\newtheorem{defn}[thm]{Definition}
\newtheorem{rem}[thm]{Remark}
\numberwithin{equation}{section}
\newcommand{\R}{\mathbb R}
\newcommand\C{\mathcal{C}}
\newcommand{\DD}{\mathbb{D}}

\newcommand{\CF}{\mathbb{F}} %
\newcommand{\dy}{\,{\rm d}y}
\newcommand{\dx}{\,{\rm d}x}

\newcommand{\dz}{\,{\rm d}z}
\newcommand{\dt}{\,{\rm d}t}
\newcommand{\ds}{\,{\rm d}s}

\newcommand{\dxi}{\,{\rm d}\xi}
\renewcommand{\div}{\mathop{\rm div}\nolimits}
\newcommand{\curl}{\mathop{\rm curl}\nolimits}
\newcommand{\om}{\omega}

\renewcommand{\th}{\theta}

\newcommand{\lec}{\lesssim}
\newcommand{\norm}[1]{\left\|#1\right\|}
\newcommand{\al}{\alpha}

\newcommand{\pd}{\partial}
\newcommand{\Om}{\Omega}
\newcommand{\ep}{\varepsilon}
\newcommand{\de}{\delta}\newcommand{\De}{\Delta}
\newcommand{\tsum}{{\textstyle \sum}}
\newcommand{\bke}[1]{\left ( #1 \right )}
\newcommand{\bkt}[1]{\left [ #1 \right ]}
\newcommand{\bket}[1]{\left \{ #1 \right \}}
\newcommand{\EQ}[1]{\begin{equation} #1 \end{equation}}
\newcommand{\EQS}[1]{\begin{equation}\begin{split} #1 \end{split}\end{equation}}
\newcommand{\EQN}[1]{\begin{equation*}\begin{split} #1 \end{split}\end{equation*}}
\newcommand{\EN}[1]{\begin{enumerate} #1 \end{enumerate}}
\newcommand{\obo}{\overline{\bm{\omega}}}

\renewcommand\rightmark%
{Gradient Estimates with the Navier Boundary Condition}

\begin{document}

		\title{Gradient estimates for the non-stationary Stokes system with the Navier boundary condition}%
		\author[H. Chen]{Hui Chen}%
		\address[H. Chen]
		{School of Science, Zhejiang University of Science and Technology, Hangzhou, 310023, People's Republic of China }
		\email{chenhui@zust.edu.cn}
		\author[S. Liang]{Su Liang}
		\address[S. Liang]
		{Department of Mathematics, University of British Columbia, Vancouver, BC V6T1Z2, Canada }
		\email{liangsu96@math.ubc.ca}
		\author[T.-P. Tsai]{Tai-Peng Tsai}
		\address[T.-P. Tsai]
		{Department of Mathematics, University of British Columbia, Vancouver, BC V6T1Z2, Canada }
		\email{ttsai@math.ubc.ca}

\maketitle

\centerline{\bf Dedicated to Professor Vladim\'{\i}r \v {S}ver\'{a}k on the occasion of his 65th birthday}

\begin{abstract}
For the non-stationary Stokes system, it is well-known that one can improve spatial regularity in the interior, but not near the boundary if it is coupled with the no-slip boundary condition. In this note we show that, to the contrary, spatial regularity can be improved near a flat boundary if it is coupled with the Navier boundary condition, with either infinite or finite slip length. The case with finite slip length is more difficult than the case with infinite slip length.

\vskip 0.2cm

		
\vskip 0.2cm
		

\end{abstract}

\section{Introduction}
Let $\Omega$ be an open subset of $\R^3$.
We consider the non-stationary Stokes equations in $\Omega\times(0,\infty)$,
\begin{align}\label{Stokes-Eqn}
	\partial_{t}\bm{u}-\Delta\bm{u}+\nabla \pi=\bm{f}+\nabla\cdot \CF,\quad \nabla\cdot \bm{u}=0,
\end{align}
where $\bm{u}=\left(u_{1},u_{2},u_{3}\right)$ and $\pi$ are the velocity field and the pressure of the fluid, respectively, and $\bm{f}=\left(f_{1},f_{2},f_{3}\right)$ and matrix $\CF=\left(F_{ij}\right)_{1\leq i,j\leq3}$ are the external forces. We use the convention that $(\nabla u)_{ij} = \pd_j u_i$ and
\begin{equation}
\label{divCF}
\nabla\cdot \CF = \div \CF, \quad (\div \CF)_i ={\textstyle \sum_j} \pd_j \CF_{ij},
\end{equation}
i.e., the derivative hits the second index. Hence $\Delta u = \div (\nabla u)$ and $u\cdot \nabla w = \div (w \otimes u)$.

When the boundary $\pd \Om$ of $\Om$ is nonempty, one needs to couple the Stokes equations \eqref{Stokes-Eqn} with a boundary condition (BC). When $\pd \Om$ is not moving, one usually assumes that the fluid does not go through the boundary, i.e.,
\begin{equation}\label{0normalBC}
	\bm{u}\cdot\bm{n}=0 ,
\end{equation}
where $\bm{n}$ is the unit outward normal vector to the boundary $\partial \Omega$. One also needs to assume certain conditions on $\bm{u}_\tau$, the tangential part of the velocity $\bm{u}$ on the boundary. There is a vast amount of literature assuming the \emph{no-slip BC}
\begin{align}\label{Dirichlet}
\bm{u} |_{\partial \Omega\times(0,\infty)} = 0,
\end{align}
which was indeed supported by most experiments before 1950s (see \cite{Neto2005}). As the experimental accuracy improved over the years, there appeared that, for certain physical situations, the no-slip BC is an accurate approximation only at macroscopic length scales; see e.g.~\cite{Lauga2007,Neto2005,Ou2004,Sochi2011}. Mathematically, it has been shown that the no-slip BC leads to the absence of collisions of rigid bodies immersed in viscous fluids \cite{Hillairet2007,MR2481302}. Navier in 1823 first proposed a slip boundary condition %
\begin{align}\label{NavierBC-general}
	\bm{u}\cdot\bm{n}=0 \quad\text{and}\quad [\left(2\DD\bm{u}+\CF\right)\bm{n}]_{\tau}+\alpha\bm{u}_{\tau}=0,
\end{align}
where $\DD\bm{u}=\frac{1}{2}\left(\nabla\bm{u}+(\nabla\bm{u})^{T}\right)$ is the \emph{deformation tensor}, and the subscript $\tau$ denotes the tangential component of a vector, i.e., $\bm{v}_{\tau}=\bm{v}-\left(\bm{v}\cdot\bm{n}\right)\bm{n}$. The constant $\alpha\geq0$ is called the \emph{friction coefficient} 
in \cite{AcevedoTapia2021,Amrouche2014}. 
The slip BC can be written as 
\[
\bm{u}_{\tau} = -\frac1\alpha\, [\left(2\DD\bm{u}+\CF\right)\bm{n}]_{\tau}
\]
which means that the tangential fluid velocity on the boundary is proportional to the reaction force of the shear stress the fluid acts upon the boundary. Hence physically $\alpha\ge 0$, which is also mathematically needed for the a priori bound.
When $\al>0$, its inverse $1/\alpha$ has the unit of length, and is called the \emph{slip length}, an important physical parameter that is measured in experiments. See \cite{JM,JBV} for derivations of the Navier condition \eqref{NavierBC-general}
as the effective boundary condition, in the limit as the length scale of the rough boundary becomes small. 
See \cite{Bocquet-Barrat,Lauga2007} for a further description of experimental and theoretical validations of the Navier BC.

When $\Omega$ is the half space $\R^3_+$, the Navier boundary condition \eqref{NavierBC-general} is reduced to
\begin{align}\label{NavierBC}
	\partial_{3}u_{1}+\CF_{13}-\alpha u_{1}=\partial_{3}u_{2}+\CF_{23}-\alpha u_{2}= u_{3}=0.
\end{align}
If in addition $\CF=-\bm{u}\otimes\bm{u}$, we recover the non-stationary Navier-Stokes equations \begin{align}\label{NS-Eqn}
	\partial_{t}\bm{u}+\bm{u}\cdot\nabla\bm{u}-\Delta\bm{u}+\nabla \pi=\bm{f},\quad \nabla\cdot \bm{u}=0
\end{align}
in $\R_{+}^3\times(0,\infty)$,
with the Navier boundary condition
\begin{align}\label{NavierBC2}
	\partial_{3}u_{1}-\alpha u_{1}=\partial_{3}u_{2}-\alpha u_{2}= u_{3}=0.
\end{align}

Many significant analytic results have been achieved for the Stokes system \eqref{Stokes-Eqn} under the Navier BC \eqref{NavierBC-general}, since the pioneering work of Solonnikov and \v{S}\v{c}adilov \cite{Solonnikov1973} on stationary boundary value problem when $\alpha=0$. We refer to \cite{AcevedoTapia2021,Amrouche2014,BeiraodaVeiga2005} for stationary Stokes and Navier-Stokes equations, \cite{AlBaba2016,Chen2010,Kucera2017} for initial-boundary value problem of the Stokes and Navier-Stokes equations with Navier BC, 
\cite{MR2325694,MR3473580} for resolvent estimates, and \cite{MR4023309} for maximal $L^p$-$L^q$ regularity of the Stokes equations under Navier BC.

In this note, we are interested in the local regularity of the solution to the Stokes equations \eqref{Stokes-Eqn} up to the flat boundary. For $x=(x_{1},x_{2},x_{3}) \in \R^{3}$, let $x^{\prime}=(x_{1},x_{2})$ be the horizontal variable and $x^{*}=(x^{\prime},-x_{3})$ the reflection point.
For $0<R$, let $B^{\prime}(R)=\bket{x'\in \R^2: |x'|<R}$, the cylinders $\C(R)=B^{\prime}(R)\times(-R,R)$ and $\C_{+}(R)=B^{\prime}(R)\times(0,R)$, and spacetime regions
  $Q(R)=\C(R)\times(-R^2,0)$ and $Q_{+}(R)=\C_{+}(R)\times(-R^2,0)$.
Denote
\begin{equation} \label{Ga.def}
\Gamma=\{x:\, |x^{\prime}|<1,\,x_{3}=0\}\subset\partial\R_{+}^3, \quad \Sigma=\Gamma\times(-1,0)\subset\partial Q_+(1).
\end{equation}

We recall some results in the literature concerning local analysis in $Q_{+}(1)$ when the boundary condition is assigned only on  $\Sigma$. Let $v(x,t)$ be a function satisfying the heat equation
and the Dirichlet BC 
\begin{equation}
\partial_{t} v-\Delta v=0 \quad \text{in}\ \ Q_{+}(1),\quad
	v\,|_{\Sigma}=0.
\end{equation}
The following a priori estimate is valid:%
\begin{equation}\label{Caccioppoli2}
	\|D_{x}^{\sigma}\, D_{t}^{m}\, v\|_{L^{2}(Q_{+}(1/2))}\leq C \|v\|_{L^{2}(Q_{+}(1))},
\end{equation}
where $D_{x}^{\sigma}$ and $D_{t}^{m}$ are usual partial derivatives with nonnegative multi-index $\sigma$ and  integer $m$. Typically, for classical linear parabolic equations, we have such smoothing estimates in space and time variables both in the interior and near a flat boundary if the coefficients of the system are sufficiently regular.

The behavior of the nonstationary Stokes system \eqref{Stokes-Eqn} is different from the heat equation.
We first consider \eqref{Stokes-Eqn} in the interior with zero external forces $f=\CF=0$ for simplicity.
By localizing maximal regularity estimates in the whole space, one can show interior estimates assuming $\bm{u},\nabla\bm{u},\pi \in L^{q,r}(Q(2))$, see e.g.~\cite[Proposition 1.1]{Seregin2011a}, which was well known before the paper. A bootstrap argument then shows that $\bm{u}$ has spatial derivatives of any order and is H\"{o}lder continuous in space-time. An alternative estimate based on the vorticity equation is given in \cite[Lemma A.2]{Chen2008a} with no assumption on the pressure term. It also improves the spatial regularity of $\bm{u}$ to arbitrary order by bootstraping. However, one still needs to assume $\pi \in L^{q,r}(Q(2))$, $r>1$, to get H\"{o}lder continuity of the velocity in time. The lack of infinite smoothing in time is shown by Serrin's example: 
$\bm{u}(x, t)=c(t) \nabla h(x)$ and $\pi(x, t)=-c^{\prime}(t) h(x)$, where $h$ is a harmonic function.

For the nonstationary Stokes system \eqref{Stokes-Eqn} in the half space under the no-slip BC \eqref{Dirichlet}, Seregin \cite[Proposition 2]{Seregin2003} showed spatial smoothing if one assumes $\pd_t \bm{u}, \nabla^2 \bm{u}, \nabla \pi \in L^{q,r}$ to start with. However, without any assumption on the pressure, the smoothing in spatial variables fails due to non-local effect of the pressure. The first counter example constructed by Kang \cite{Kang2004a} shows that the gradient of a continuous weak solution may be unbounded. A simplified example was constructed in Seregin-\v Sver\'ak \cite{Seregin2011a} as a shear flow.
Chang and Kang \cite[Theorem 1.1]{Chang2021} constructed a bounded very weak solution whose derivatives are unbounded in any $L^{p}$ with $1<p<\infty$. See \cite{Chang2020,KLLT1,KLLT2,Chang2023,KangMin} for related study along this line  and \cite{Gustafson2006,Seregin2002} for partial regularity results on the boundary with no-slip BC.

The goal of this paper is to show that, for the nonstationary Stokes system \eqref{Stokes-Eqn} in the half space under the Navier BC \eqref{NavierBC}, even a very weak solution has certain smoothing effect in the spatial variable. This is completely opposite to the case with no-slip BC \eqref{Dirichlet}. See next section for exact definitions of very weak solutions and weak solutions. Similar study was initiated by Dong, Kim and Phan \cite{DKP}, where boundary second derivative estimates for generalized Stokes systems with VMO coefficients under Navier BC with $\alpha=0$ were proved. For its relevance to our results, see Comment (iii) on Theorem \ref{thm2} and Comment (i) on Theorem \ref{thm0}.

\begin{thm}[Boundary gradient estimate]\label{thm1} 
Assume that $\bm{u}$ is a \emph{very weak solution} to the Stokes equations \eqref{Stokes-Eqn} in $Q_{+}(1)$ with Navier boundary condition \eqref{NavierBC} for a constant $\alpha\geq 0$ on the flat boundary $\Sigma$. Suppose that $\bm{u},\CF\in L^{q,r}\left(Q_{+}(1)\right)$ and $\bm{f}\in L^{q_{*},r}\left(Q_{+}(1)\right)$ for $1<q,r<\infty$ and $q_{*}=\left\{\begin{array}{lll}
	\max\{1,\frac{3q}{q+3}\},&if\, q\neq \frac{3}{2}\\
	1+,&if\,q=\frac{3}{2}
\end{array}\right.$. Then $\bm{u}\in L^{r}\left(-\frac{1}{4},0;W^{1,q}(\C_{+}(\frac{1}{2}))\right)$ is a \emph{weak solution}, and there is a constant $C$ depending only on $q,r,\alpha$ such that
\begin{equation}\label{main1}
\| \nabla\bm{u}\|_{L^{q,r}\left(Q_{+}(1/2)\right)}\leq C\|\bm{u}\|_{L^{q,r}\left(Q_{+}(1)\right)}+C\|\bm{f}\|_{L^{q_{*},r}\left(Q_{+}(1)\right)}+C\|\CF\|_{L^{q,r}\left(Q_{+}(1)\right)}.
\end{equation}
\end{thm}

\emph{Comments on Theorem \ref{thm1}:}
\begin{enumerate}[(i)]
\item There is no assumption on the pressure $\pi$ in Theorem \ref{thm1}. It is similar to the interior estimate in \cite[Lemma A.2]{Chen2008a}, and the stationary cases in \cite{MR1789922,Kang2004}. %

\item  For most analytic results in the literature, what is true for no-slip BC is usually also true for Navier BC. However,
estimate \eqref{main1} is wrong if the Navier BC \eqref{NavierBC} is replaced by the no-slip BC \eqref{Dirichlet}, see \cite{Kang2004a,Seregin2011a,Chang2020,Chang2021,KLLT2}. This is surprising and the main merit of this paper.

\item Another property that depends on the boundary condition is the  occurrence of finite-time collision of solids inside viscous incompressible fluids. For example, in the simple geometry that a spherical solid moves toward a flat boundary, there is no finite-time collision under no-slip BC \cite{Hillairet2007,MR2481302}, but there is  finite-time collision under Navier BC \cite{MR3281946}.\footnote{TT learned this difference from a talk of 
Prof.~\v{S}\'{a}rka Ne\v{c}asov\'{a} in 2017, which inspired this project.} This problem is also considered under different regularity of the boundary \cite{MR2592281}, and under mixed BC in a bounded Lipschitz domain \cite{MR3567970}. For the compressible fluid case, see 
the references in \cite{JinNecasova}.

\item A key difficulty under the no-slip BC is that both the vorticity and the pressure have no boundary condition, although they satisfy the nonhomogeneous heat equation and the Poisson equation, respectively. When it comes to Navier BC, however, it is well known that the vorticity satisfies a boundary condition for both $\al =0$ and $\al>0$, see 
Chen-Qian \cite[(1.4)]{Chen2010}. When $\Om=\R^3_+$ with a flat boundary and assume $\CF=0$ for simplicity, a vorticity BC can be derived easily from \eqref{NavierBC}:
\begin{equation}\label{NavierBComega}
\om_{1}= -\al u_2, \quad
\om_{2}= \al u_1, \quad
\partial_{3}\om_{3}=\alpha \om_{3}.
\end{equation}

\item The existence of a vorticity BC under the Navier BC is a key difference with the no-slip BC, and is the heuristic reason why \eqref{main1} may be correct. We will give a proof of Theorem \ref{thm1} utilizing the 
vorticity BC \eqref{NavierBComega} in the Appendix, assuming $\bm{u}$ has sufficient regularity to make sense of the heat equation and the BC. However, in our main proof in Section \ref{Sec 4}, we will only assume $\bm{u}$ is a very weak solution and will not use \eqref{NavierBComega}. The main proof is based on choosing suitable test fields $\Phi$ in the definition of very weak solutions to get \emph{uniform gradient bound} of approximations of $\bm{u}$. 

\end{enumerate}

In the next theorem we show estimates of higher order derivatives, assuming the integrability of the pressure. \begin{thm}[Higher-order derivative estimates]\label{thm2}
Let $\bm{f},\CF,\div\CF\in L^{q,r}\left(Q_{+}(1)\right)$, $1<r,q<\infty$, with $\CF_{13}=\CF_{23}=0$ on the boundary $\Sigma$.
If $\left(\bm{u},\pi\right)\in L^{r}\left(-1,0;W^{1,q}(\C_{+}(1))\right)\times L^{q,r}\left(Q_{+}(1)\right)$ is a weak solution pair to the Stokes equations \eqref{Stokes-Eqn} in $Q_{+}(1)$ with Navier boundary condition \eqref{NavierBC} for a constant $\alpha \geq 0$ on the flat boundary $\Sigma$, then $\partial_{t}\bm{u},\nabla^2\bm{u}, \nabla\pi\in L^{q,r}\left(Q_{+}(1/2)\right)$ and there is a constant $C$ depending only on $q,r,\alpha$ such that
\begin{align}\label{main2}
	&\| \partial_{t}\bm{u}\|_{L^{q,r}\left(Q_{+}(1/2))\right)}+	\| \nabla^2\bm{u}\|_{L^{q,r}\left(Q_{+}(1/2))\right)}+\|\nabla\pi\|_{L^{q,r}(Q_{+}(1/2))}
	\notag\\
\leq& C\|\bm{u}\|_{L^{q,r}(Q_{+}(1))}+C\|\pi\|_{L^{q,r}(Q_{+}(1))}+C\|\bm{f}+\div\CF\|_{L^{q,r}(Q_{+}(1))}.
\end{align}
If in addition $\nabla \left(f+\div\CF\right)\in  L^{q,r}\left(Q_{+}(1)\right)$, we have
 \begin{align}\label{main3}
	&\|\nabla\partial_{t}\bm{u}\|_{L^{q,r}(Q_{+}(1/2))}+\|\nabla^3 \bm{u}\|_{L^{q,r}(Q_{+}(1/2))} +\|\nabla^2\pi\|_{L^{q,r}(Q_{+}(1/2))}\notag\\
	\lesssim&\|\bm{u}\|_{L^{q,r}(Q_{+}(1))}+\|\pi\|_{L^{q,r}(Q_{+}(1))}
	+\|\bm{f}+\div\CF\|_{L^{r}\left(-1,0;W^{1,q}(\C_{+}(1))\right)}.
\end{align}
\end{thm}
\emph{Comments on Theorem \ref{thm2}:}
\begin{enumerate}[(i)]
\item In Theorem \ref{thm2} we make a key assumption that $\pi\in L^{q,r}\left(Q_{+}(1)\right)$. 
\item The assumption $\CF \in L^{q,r}$ is used to make sense of the weak form \eqref{weakform2}. We do not need $\norm{\CF}_{L^{q,r}}$ in the right sides of \eqref{main2} and \eqref{main3}.

\item When one further assumes $\nabla^2\bm{u}\in L^{q,r}\left(Q_{+}(1)\right)$ and $\partial_{t}\bm{u}, \nabla\pi\in L^{1}\left(Q_{+}(1)\right)$, one can derive the bound of $\| \nabla^2\bm{u}\|_{L^{q,r}\left(Q_{+}(1/2))\right)}$ in \eqref{main2} by applying \cite[Theorem 1.2]{DKP} to $\bm{v} = e ^{-\alpha x_3} \bm{u}$ and $p = e ^{-\alpha x_3}\pi$, noting that \cite[Theorem 1.2]{DKP}  allows nonzero divergence as a source term. 
As this change of variables introduces the term $ \| \pi \|_{L^{q,r}\left(Q_{+}(1))\right)}$ on the right side, it does not recover \eqref{main1}.

\end{enumerate}

Now we consider the special case $\alpha=0$. By introducing an ``even-even-odd" extension of the velocity (already used in, e.g., \cite{DKP,MR4558989}), we can obtain arbitrary spatial smoothness of the velocity as long as the eternal forces $\bm{f}$ and $\CF$ satisfy enough regularity and boundary conditions on $\Sigma$. In order to better illustrate this point, we consider a simple case $f=0$ and $\CF=-\bm{u}\otimes\bm{u}$. 

\begin{thm}[Infinite slip length case]\label{thm0}
	Assume that $\bm{u}$ is a very weak solution to the Navier-Stokes equations \eqref{NS-Eqn} in $Q_{+}(1)$ with $\bm{f}=0$ and Navier boundary condition \eqref{NavierBC2} with $\alpha=0$ on the flat boundary $\Sigma$. If $\bm{u}\in L^{\infty}(Q_{+}(1))$, then $\nabla^n\bm{u}\in L^{\infty}(Q(1/4))$ for all $n\geq1$.
\end{thm}

\emph{Comments on Theorem \ref{thm0}:}
\begin{enumerate}[(i)]
\item By Theorem \ref{thm1} and Lemma \ref{lem3}, the ``even-even-odd" extension of the solution $\bm{u}$ in Theorem \ref{thm0}
is a weak solution in $Q(1/2)$, and the theorem follows from the usual interior regularity theory.
Although this extension when $\alpha=0$ is well known to experts, e.g., authors of \cite{DKP,MR4558989}, it takes some care to deal with very weak solutions.
See \cite{Luo2015} for more Serrin-type interior regularity criteria for very weak solutions of \eqref{NS-Eqn}.

\end{enumerate}

Our paper is organized as follows. We recall some notations and preliminary results in Section \ref{Sec 2}. We prove Theorem \ref{thm1} in Section \ref{Sec 4}, Theorem \ref{thm2} in Section \ref{Sec 5}, and Theorem \ref{thm0} in Section \ref{Sec 3}. An alternative proof of Theorem \ref{thm1} assuming more regularity is in Appendix \S\ref
{sec:app}.

\section{Notations and preliminaries}\label{Sec 2}
We first specify the notations.
For two comparable quantities, the inequality $X\lesssim Y$ stands for $X\leqslant C Y$ for some positive constant $C$. The dependence of the constant $C$ on other parameters or constants are usually clear from the context, and we will often suppress this dependence.
We have already defined $x'$, $x^*$, $B^{\prime}(R)$, $\C(R)$, $\C_{+}(R)$, $Q(R)$,  $Q_{+}(R)$, $\Gamma$ and $\Sigma$ in the paragraph above \eqref{Ga.def}.

For an open domain $\Omega\subset \R^3$, $1\leq q,r \leq \infty$ and integer $m\geq 1$, let $L^{q}(\Omega)$ and $W^{m,q}(\Omega)$ be the usual (vector-valued) Lebesgue and Sobolev spaces, respectively. Denote by $L^{q,r}\left(\Omega\times I\right)=L^{r}\left(I;L^{q}(\Omega)\right)$ a mixed Lebesgue space for $I\subset \R$ and  $L^{q}\left(\Omega\times I\right)=L^{q,q}\left(\Omega\times I\right)$. Let $C_{c}^{2,1}\left(\Omega\times I\right)=\{g\in C_{c}\left(\Omega\times I\right)|\nabla^2 g, \,\partial_{t}g\in C\left(\Omega\times I\right)\}$ and $C_{c}^{1}\left(\Omega\times I\right)=\{g\in C_{c}\left(\Omega\times I\right)|\nabla g, \,\partial_{t}g\in C\left(\Omega\times I\right)\}$.%

We write $D_{x_{i}}=\frac{\partial}{\partial x_{i}}$ the partial derivative  with respect to the variable $x_{i}, 1 \leq i \leq 3$, $D_{x}^{\sigma}=D_{x_{1}}^{\sigma_{1}}D_{x_{2}}^{\sigma_{2}}D_{x_{3}}^{\sigma_{3}}$ for a positive multi-index $\sigma=\left(\sigma_{1},\sigma_{2},\sigma_{3}\right)$ with $|\sigma|=\sigma_{1}+\sigma_{2}+\sigma_{3}$ and the gradient $\nabla_{x}=\left(D_{x_{1}},D_{x_{2}},D_{x_{3}}\right)$.

We now introduce some cut-off functions which will be used throughout this paper. Fix a cut-off function $\phi\in C_{c}^{\infty}[0,3/4)$ with $\phi=1$ on $[0,1/2]$, and let $\zeta(x)=\phi(|x^{\prime}|)\phi(|x_{3}|)\in C_{c}^{\infty}(\C(1))$, and $\theta(t)\in C_{c}^{\infty}(-1,0]$ satisfying $\theta(t)=1$ for $t\in(-1/4,0)$.%

 For $0<\varepsilon<\frac{1}{100}$, denote $\eta_{\varepsilon}(x)=\varepsilon^{-3}\eta(\frac{x}{\varepsilon})$ where $\eta(x)=\eta(|x|)$ is a cut-off function satisfing $\int_{\R^3}\eta(x)\dx=1$, and%
\begin{equation}
	\eta_{\ep}^{+}(x,y)=\eta_{\ep}(x-y)+\eta_{\ep}(x-y^*),\quad \eta_{\ep}^{-}(x,y)=\eta_{\ep}(x-y)-\eta_{\ep}(x-y^*).
\end{equation}
For $g(x)\in L^{1}_{\text{loc}}(\overline{\R^3_{+}})$, define the even and odd extensions of $g$ to $\R^3$, 
\begin{equation}\label{Eg3}
	g^+(x)=\left\{\begin{array}{lll}
		g(x),&x_{3}>0\\
		g(x^*),&x_{3}<0
	\end{array}\right.,\qquad	g^-(x)=\left\{\begin{array}{lll}
		g(x),&x_{3}>0\\
		-g(x^*),&x_{3}<0
	\end{array}\right.,
\end{equation}
and the operators defined for $x\in\R^3$ for
their mollifications,
\begin{align}\label{Eg+}
	E^{+}(g)(x)&=\int_{\R^3}g^+(x-y)\eta_{\ep}(y)\dy=\int_{\R^3_{+}}g(y)\eta_{\ep}^{+}(x,y)\dy,\\
E^{-}(g)(x)&=\int_{\R^3}g^-(x-y)\eta_{\ep}(y)\dy=\int_{\R^3_{+}}g(y)\eta_{\ep}^{-}(x,y)\dy.
\label{Eg-}
\end{align}

The following are some elementary properties. 
We skip their proof.

\begin{lem}\label{lem0} We have
	\begin{enumerate}
		\item[(i)] $D_{x_{3}}^{k}\zeta(x^{\prime},0)=0$ for any integer $k\geq 1$;
		\item[(ii)] $D_{x_{k}}\eta^{\pm}_{\ep}(x,y)=-D_{y_{k}}\eta_{\ep}^{\pm}(x,y)$ for $k=1,2$ and $D_{x_{3}}\eta_{\ep}^{\pm}(x,y)=-D_{y_{3}}\eta_{\ep}^{\mp}(x,y)$,\\
		$($a shorter form: $D_{x_{k}}\eta^{\pm}_{\ep}=-D_{y_{k}}\eta_{\ep}^{\epsilon_k\pm}$, $\epsilon_k=1-2\de_{k3})$;
		\item[(iii)] $E^{+}(g),E^{-}(g)$ are smooth in $\R^3$ and converge to $g(x)$ in $L^{q}(\R_{+}^3)$  as $\ep \to 0$ for  $1\leq q <\infty$;
		\item[(iv)] $\|E^{+}(g)\|_{L^{q}(\R^3)}+\|E^{-}(g)\|_{L^{q}(\R^3)}\lesssim\|g\|_{L^{q}(\R^3_{+})}$;
	\end{enumerate}
\end{lem}

\subsection{Definitions of weak and very weak solutions}
For divergence free vector function $\bm{u}\in L^{q}(\C_{+}(1))$, its normal boundary value
$\bm{u}\cdot\bm{n}$ is well-defined. Here $\bm{n}$ denotes the unit outer normal vector of the boundary of $\C_+(1)$.
We say $\bm{u}\cdot\bm{n}=0$ on the flat boundary $\Gamma$, if for each scalar function $\Psi \in C_{c}^{1}\left(\C_{+}(1)\cup \Gamma\right)$, which may not vanish on the boundary $\Gamma$,
\begin{align}\label{normalBC}
	\int_{\mathbb{R}_{+}^{3}} \bm{u} \cdot \nabla \Psi\dx=0.
\end{align}

Next, we introduce notions of very weak solutions and weak solutions for the non-stationary Stokes equations with Navier boundary condition on the flat boundary.
See \cite[Definition 3.10, p.269]{AcevedoTapia2021} for weak solutions of the stationary Stokes equations with Navier BC on general boundary. Let's motivate the weak formulation: Recall in \eqref{Stokes-Eqn}, $f$ denotes a vector field and $\CF$ a matrix function. Also recall the convention \eqref{divCF}. Using $\div \bm{u}=0$, we have $\Delta \bm{u} = 2 \div \DD \bm{u}$ and can rewrite \eqref{Stokes-Eqn} as
\[
\partial_{t}\bm{u}- \div(2 \DD \bm{u} +\CF)+\nabla \pi=\bm{f},\quad \nabla\cdot \bm{u}=0.
\]
Testing it with any divergence free test field $\Phi(x,t)$ that vanishes at $t=0,T$ and  $\Phi \cdot \bm{n}=0$ on $\pd \Omega\times (0,T)$, we get
\begin{multline}\label{weakform0}
\int_0^T \!\int_\Om\bket{ -\bm{u} \cdot \partial_{t}\Phi+\bke{2\DD \bm{u}+\CF}:\nabla \Phi -  f\cdot \Phi } \dx\dt
\\
=\int_0^T \!\int_{\pd\Om} \bke{2\DD \bm{u}+\CF}\bm{n}\cdot \Phi \,dx^{\prime}\dt
=- \int_0^T \!\int_{\pd\Om} \alpha \bm{u}\cdot \Phi \,dx^{\prime}\dt.
\end{multline}
In the last equality we have used the Navier BC \eqref{NavierBC}. The equality \eqref{weakform0} can be used to define weak solutions with Navier BC in general domains. There is an alternative form
\begin{equation}\label{05160}
\int_\Omega 2\DD \bm{u}:\nabla \Phi = 2\int_\Omega \DD \bm{u}:\DD\Phi ,
\end{equation}
by switching indexes. It is convenient for getting the a priori bound for $\int_\Omega |\DD \bm{u}|^2$ with $\Phi=\bm{u}$.

When $\Om=\R^3_+$, there is another form
\begin{equation}\label{0516a}
\int_\Omega 2\DD \bm{u}:\nabla \Phi = \int_\Omega \nabla \bm{u}:\nabla \Phi .
\end{equation}
One advantage is that it gives the a priori bound for the full gradient $\int_\Omega |\nabla \bm{u}|^2$ with $\Phi=\bm{u}$.
Indeed, the difference of the two sides of \eqref{0516a}
\[
J
=\int_\Om (\pd_i u_j )\pd_j \Phi_i = \int_\Om \pd_j (\Phi_i \pd_i u_j )= \int_{\pd\Om} n_j (\Phi_i \pd_i u_j  )= \int_{\pd\Om} \bm{n} \cdot (\Phi \cdot \nabla \bm{u}).
\]
When $\Om=\R^3_+$, 
\[
J = -\int_\Gamma \Phi_i \pd_i u_3 = - \int_\Gamma \Phi_1 \pd_2 u_3+\Phi_1 \pd_2 u_3+ \Phi_3 \pd_3 u_3 = 0,
\]
as $\pd_2 u_3= \pd_2 u_3= \Phi_3 =0$ on $\Gamma$. However, \eqref{0516a} is not true for general $\Om$ with a curved boundary. For example, let $\Om=B_1$
 in $\R^3$ with unit outer normal $\bm{n}=re_r+ze_z$ on $\pd \Om$ in cylindrical coordinates $r,\theta,z$. Let $\bm{u}=\Phi=f(r,z)e_\theta$. Then $\div \bm{u}=0$ and $\bm{u}\cdot \bm{n}=0$ on $\pd \Om$. We have
\[
\bm{u} \cdot \nabla \bm{u} = f \frac 1r \pd_\theta (f e_\theta) = \frac {f^2}r \pd_\theta  e_\theta = -  \frac {f^2}r e_r.
\]
Thus
\[
\int_{\pd\Om} \bm{n} \cdot (\bm{u} \cdot \nabla \bm{u}) = \int_{\pd\Om} (re_r+ze_z) \cdot (-  \frac {f^2}r e_r) = - \int_{\pd\Om} f^2
\]
which is negative if $f$ does not vanish entirely on $\pd B_1$.

We now define very weak solutions and weak solutions in $Q_{+}(1) \subset \R^3_+ \times (-1,0)$.

\begin{defn}\label{def:vws}
Let $\bm{f},\CF\in L^{1}\left(Q_{+}(1)\right)$. A vector field $\bm{u}$ is called a \emph{very weak solution} to the Stokes equations \eqref{Stokes-Eqn} in $Q_{+}(1)$ with Navier boundary condition  \eqref{NavierBC} on the flat boundary $\Sigma$, if
\begin{enumerate}
\item[(i)] $\bm{u} \in L^{1}\left(Q_{+}(1)\right)$  satisfies \eqref{normalBC} for a.e.~$t$, (hence $\bm{u}$ is divergence free and $\bm{u}\cdot \bm{n}=0$ on $\Sigma$);

\item[(ii)]  For each divergence free vector $\Phi=\left(\Phi_{1},\Phi_{2},\Phi_{3}\right) \in C_{c}^{2,1}\left(Q_{+}(1) \cup \Sigma\right)$ with
\begin{align}\label{NavierBC3}
	\partial_{3}\Phi_{1}-\alpha \Phi_{1}= \partial_{3}\Phi_{2}-\alpha \Phi_{2}=
	\Phi_{3}=0
\end{align}
on the boundary $\Sigma$, we have
\begin{align}\label{Stokes-weak-form}
	- \iint_{Q_{+}(1)} \bm{u} \cdot \left(\partial_{t}+\Delta\right) \Phi \dx\dt= \iint_{Q_{+}(1)} \bm{f}\cdot\Phi-\CF:\nabla \Phi \dx\dt.
\end{align}
\end{enumerate}
\end{defn}

\begin{defn}\label{def:ws}
Let $\bm{f},\CF\in L^{1}\left(Q_{+}(1)\right)$. A vector field $\bm{u}$ is called a \emph{weak solution} to the Stokes equations \eqref{Stokes-Eqn} in $Q_{+}(1)$ with Navier boundary condition \eqref{NavierBC} on the flat boundary $\Sigma$, if
\begin{enumerate}
	\item[(i)] $\bm{u} \in L^{1}\left(-1, 0; W^{1,1}(\C_{+}(1))\right)$ is divergence free with $u_{3}|_{\Sigma}=0$;
	
	\item[(ii)]  for each divergence free vector $\Phi=\left(\Phi_{1},\Phi_{2},\Phi_{3}\right) \in C_{c}^{1}\left(Q_{+}(1) \cup \Sigma\right)$ with $\Phi_{3}(x^{\prime},0,t)=0$, we have
\begin{equation}\label{weakform}
\begin{split}
	&\quad\ \iint_{Q_{+}(1)} -\bm{u} \cdot \partial_{t}\Phi+\nabla \bm{u}:\nabla \Phi \dx\dt+\tsum_{k=1}^2\iint_{\Sigma}\alpha u_{k}\cdot\Phi_{k}\dx^{\prime}\dt\\
	&=\iint_{Q_{+}(1)} \bm{f}\cdot\Phi-\CF:\nabla\Phi \dx\dt.
\end{split}
\end{equation}
\end{enumerate}
\end{defn}

\begin{defn}\label{def:wsp}
Let $\bm{f},\CF\in L^{1}\left(Q_{+}(1)\right)$. A pair $\left(\bm{u},\pi\right)$ is called a \emph{weak solution pair} to the Stokes equations \eqref{Stokes-Eqn} in $Q_{+}(1)$ with Navier boundary condition \eqref{NavierBC} on the flat boundary $\Sigma$, if
\begin{enumerate}
	\item[(i)] $\bm{u} \in L^{1}\left(-1, 0; W^{1,1}(\C_{+}(1))\right)$ is divergence free with $u_{3}|_{\Sigma}=0$, $\pi\in L^{1}\left(Q_{+}(1)\right)$;
	
	\item[(ii)]  for each vector $\Phi=\left(\Phi_{1},\Phi_{2},\Phi_{3}\right) \in C_{c}^{1}\left(Q_{+}(1) \cup \Sigma\right)$ with $\Phi_{3}(x^{\prime},0,t)=0$, we have
	\begin{equation}\label{weakform2}
		\begin{split}
			&\quad \ \iint_{Q_{+}(1)} -\bm{u} \cdot \partial_{t}\Phi+\nabla \bm{u}: \nabla \Phi -\pi\cdot \mathrm{div}\Phi\dx\dt+\tsum_{k=1}^2\iint_{\Sigma}\alpha u_{k}\cdot\Phi_{k}\dx^{\prime}\dt\\
			&=\iint_{Q_{+}(1)} \bm{f}\cdot\Phi-\CF:\nabla\Phi \dx\dt.
		\end{split}
	\end{equation}
\end{enumerate}
\end{defn}

Note the different classes of test fields $\Phi$ in these definitions: For very weak solutions, $\Phi$ is divergence-free with Navier BC. For weak solutions, $\Phi$ is divergence-free with zero normal on boundary.  For weak solution pairs, $\Phi$ has zero normal on boundary, but may not be divergence-free.

\subsection{A differentiable very weak solution is a weak solution}

A weak solution is clearly a very weak solution. In this subsection we show that a very weak solution of Stokes equations with Navier BC and sufficient regularity is a weak solution. 

Recall the study of the oblique derivative boundary condition in Gilbarg and Trudinger \cite[\S6.7]{Gilbarg2001} for the Poisson equation and Keller \cite{Keller1981} for heat and wave equations.
Inspired by the Green function formulas for Robin BC, \cite[(6.62), p.121]{Gilbarg2001} for the Poisson equation and  \cite[(5.3), p.299]{Keller1981} for the heat equation, we define the operator $I_\al$: 
For $x\in \R^3_+$, let
\begin{equation}
I_{\alpha}(g)(x)=\int_{0}^{\infty}e^{-\alpha\xi}\int_{\R^3_{+}} g(y)\cdot\eta_{\ep}(x^{\prime}-y^{\prime},x_{3}+y_{3}+\xi)\dy\dxi.\label{Ig}
\end{equation}
It has the good properties that
\begin{equation} \label{eq2.14}
g^\ep = E^+ g - 2\al I_\alpha g \quad \text{
satisfies the Robin BC} \quad(D_3 - \al )g^\ep|_{x_3=0}=0,
\end{equation}
and $I_\alpha g$ is a local perturbation as it is supported in $x_3\le\ep$ and it vanishes in $L^q$ as $\ep \to 0$.
It will be used to construct suitable test fields $\Phi^\ep$ in \eqref{eq2.17} for the weak form \eqref{Stokes-weak-form}.

The following are some elementary properties of $I_\alpha$. 
We skip the proof.
\begin{lem}\label{lem0a} We have
\EN{
\item[(i)] $I_{\alpha}(g)(x)=\int_{0}^{\ep}e^{-\alpha\xi}\int_{\R^3_{+}} g(y)\cdot\eta_{\ep}(x^{\prime}-y^{\prime},x_{3}+y_{3}+\xi)\dy\dxi$
is smooth in $\R^3_{+}$ with $I_{\alpha}(g)(x)=0$ if $x_{3}>\ep$, $D_{x_{i}}I_{\alpha}(g)= I_{\alpha}(D_{x_{i}}g)$ for $i=1,2$, and
		\begin{align}\label{2.5}
					 D_{x_{3}}I_{\alpha}(g)&=\alpha I_{\alpha}(g) -\int_{\R^3_{+}}g(y)\eta_{\ep}(x-y^*)\dy
					\\
		&= -I_{\alpha}(D_{x_{3}}g) -\int_0^\infty\!\!e^{-\alpha\xi} \int_{\R^2}g(y',0)\eta_{\ep}(x^{\prime}-y^{\prime},x_{3}+\xi)\dy'd\xi. \label{2.6}
		\end{align}
	\item[(ii)] For $g\in L^{\infty}(\R_{+}^3)$, we have
	\begin{align}
		\|I_{\alpha}(g)\|_{L^{q}(\C_{+}(1))}\lesssim  \ep^{1+\frac{1}{q}} \,\|g\|_{L^{\infty}(\R_{+}^3)},
	\end{align}
    and
    \begin{align}\label{2.7}
    	\|D_{x_{3}}I_{\alpha}(g)\|_{L^{q}(\C_{+}(1))}+\|\int_{\R^3_{+}}g(y)\eta_{\ep}(x-y^*)\dy\|_{L^{q}(\C_{+}(1))}\lesssim \ep^{\frac{1}{q}}\,\|g\|_{L^{\infty}(\R_{+}^3)}.
    \end{align}
}   
\end{lem}

The following lemma shows that a very weak solution of Stokes equations is a weak solution if it has sufficient regularity.

\begin{lem}\label{lem2}
Assume that $\bm{f},\CF\in L^{q,1}(Q_{+}(1))$ for $1<q<\infty$ and $\bm{u}\in L^{1}(-1,0;W^{1,q}(\C_{+}(1)))$ is a very weak solution to the Stokes equations \eqref{Stokes-Eqn} in $Q_{+}(1)$ with Navier boundary condition \eqref{NavierBC} on the flat boundary $\Sigma$. Then $\bm{u}$ is a weak solution.
\end{lem}
\begin{proof}
Let $1<q^{\prime}<\infty$ with $\frac{1}{q}+\frac{1}{q^{\prime}}=1$ and  $\Phi=\left(\Phi_{1},\Phi_{2},\Phi_{3}\right) \in
C_{c}^{1}(Q_{+}(1) \cup \Sigma)$ be any given divergence free vector field with $\Phi_{3}(x^{\prime},0,t)=0$. The main idea is to construct a sequence of divergence free vectors $\Phi^{\ep}\in C_{c}^{2,1}\left(Q_{+}(1)\cup\Sigma\right)$, which satisfy the boundary condition \eqref{NavierBC3} and converge to $\Phi$ as $\ep\rightarrow0_+$.

Let the approximation $E(\Phi)=\left(E^+ (\Phi_1), E^+ (\Phi_2), E^- (\Phi_3)\right)$ where the operators $E^{+},E^{-}$ are defined in \eqref{Eg+} and \eqref{Eg-} for a sufficiently small $\ep>0$.
Notice that  $E(\Phi)\in C_{c}^{2,1}(Q(1))$ is divergence free with $E^{-}(\Phi_{3})|_{x_{3}=0}=0$. We seek a correction
$\Phi^{\ep}=\left(\Phi_{1}^{\ep},\Phi_{2}^{\ep},\Phi_{3}^{\ep}\right)$ with
\begin{align}
 	\Phi^{\ep}_{1}(x,t)&=E^{+}\left(\Phi_{1}\right)-2\alpha I_{\alpha}(\Phi_{1})+\Psi_{1}^{\ep},\notag\\
 	\Phi^{\ep}_{2}(x,t)&=E^{+}\left(\Phi_{2}\right)-2\alpha I_{\alpha}(\Phi_{2})+\Psi_{2}^{\ep},\label{eq2.17} \\
 	\Phi^{\ep}_{3}(x,t)&=E^{-}\left(\Phi_{3}\right)+2\alpha I_{\alpha}(\Phi_{3})+\Psi_{3}^{\ep},\notag
\end{align}
where $0\leq x_{3}<1$ and the operator $I_{\alpha}$ is defined in \eqref{Ig}. Clearly $2\al I_\al( \Phi_1,  \Phi_2,- \Phi_3) $ is divergence free.
The correction term $\Psi^{\ep}=\left(\Psi_{1}^{\ep},\Psi_{2}^{\ep},\Psi_{3}^{\ep}\right)$ should be a divergence free vector such that $\Psi^\ep_1$ and $\Psi^\ep_2$ are supported away from $\Sigma$ and 
$\Psi_{3}^{\ep}=-2\alpha I_{\alpha}(\Phi_{3})$ on $\Sigma$ so that $\Phi^\ep$ satisfy the BC \eqref{NavierBC3} using \eqref{eq2.14}.
We choose
\begin{equation}\label{2.13}
		\Psi_{3}^{\ep}(x,t)=-2\alpha I_{\alpha}(\Phi_{3})(x^{\prime},0,t)\phi(x_{3}),
\end{equation}
(so that $\Phi^{\ep}_{3}(x',0,t)=0$)
and
\begin{equation}\label{divegencefree}
	D_{x_{1}}\Psi_{1}^{\ep}+D_{x_{2}}\Psi_{2}^{\ep}=-D_{x_{3}}\Psi_{3}^{\ep}=2\alpha I_{\alpha}(\Phi_{3})(x^{\prime},0,t)D_{x_{3}}\phi(x_{3}),
\end{equation}
where $D_{x_{3}}\Psi_{3}^{\ep}\in C_{c}^{2,1}(Q_{+}(1))$.
Since $\Phi_{3}(x^{\prime},0,t)=0$ and
\begin{equation*}
	D_{x_{3}}\int_{\R^2}\Phi_{3}\dx^{\prime}=-\int_{\R^2}D_{x_{1}}\Phi_{1}+D_{x_{2}}\Phi_{2}\dx^{\prime}=0
\end{equation*}
due to divergence free property of the vector $\Phi$, we have 
\begin{align}\label{2.9}
	\int_{\R^2}\Phi_{3}\dx^{\prime}=0,
\end{align}
for all $0\leq x_{3}<1$ and $t\in(-1,0)$. Thus, letting $z'=x'-y'$ and by \eqref{2.9}, we have
\begin{align*}
&	\int_{\R^2}I_{\alpha}(\Phi_{3})(x^{\prime},0,t)\dx^{\prime}=\int_{\R^2}\int_{0}^{\infty}e^{-\alpha\xi}\int_{\R^3_{+}} \Phi_{3}(y,t)\,\eta_{\ep}(x^{\prime}-y^{\prime},y_{3}+\xi)\dy\dxi\dx^{\prime}
	\\
	&=\int_{0}^{\infty}e^{-\alpha\xi}\int_{0}^{\infty}\int_{\R^2}  \left(\int_{\R^2}\Phi_{3}(x'-z',y_3,t)\dx^{\prime}\right)\eta_{\ep}(z^{\prime},y_{3}+\xi)\dz^{\prime}\dy_{3}\dxi
	=0.
\end{align*}
Hence, $\int_{\R^2}D_{x_{3}}\Psi_{3}^{\ep}(x,t)\dx^{\prime}=0$.
In order to find $\Psi_{1}^{\ep},\Psi_{2}^{\ep} \in C_{c}^{2,1}(Q_{+}(1))$ solving the equation \eqref{divegencefree}, we introduce the 2D Bogovski\u{\i} formula
\begin{equation}\label{Bogovskii}
\left(\Psi_{1}^{\ep},\Psi_{2}^{\ep}\right)(x,t)=2\alpha\int_{\R^2}I_{\alpha}(\Phi_{3})(y^{\prime},0,t) N(x^{\prime},y^{\prime}) \dy^{\prime}\cdot D_{x_{3}}\phi(x_{3})
\end{equation}
with
\begin{align}
	N(x^{\prime}, y^{\prime})=\frac{1}{\int_{B^{\prime}(1)}\phi(|y^{\prime}|)\dy^{\prime}}\cdot\frac{x^{\prime}-y^{\prime}}{|x^{\prime}-y^{\prime}|^{2}} \int_{|x^{\prime}-y^{\prime}|}^{\infty} \phi\left(\left|y^{\prime}+s \frac{x^{\prime}-y^{\prime}}{|x^{\prime}-y^{\prime}|}\right|\right) s \ds.
\end{align}
See \cite[Lemma \textrm{III}.3.1]{Galdi2011}.

\medskip

$\bullet$ We claim that $\Phi^{\ep}\in C_{c}^{2,1}(Q_{+}(1)\cup\Sigma)$ are divergence free and satisfy the boundary condition \eqref{NavierBC3}.
By Lemma \ref{lem0}(ii) and Lemma \ref{lem0a}(i), in particular \eqref{2.6} with $\Phi_3(y',0,t)=0$) and \eqref{divegencefree}, it is easy to verify that
\begin{align*}
	\nabla \cdot \Phi^{\ep}=E^{+}\left(\nabla\cdot\Phi\right)-2\alpha I_{\alpha}(\nabla\cdot\Phi)+\nabla\cdot\Psi^{\ep}=0.
\end{align*}
For $k=1,2$,
\begin{align*}
	D_{x_{3}}\Phi_{k}^{\ep}-\alpha\Phi_{k}^{\ep}&=E^{-}\left(D_{x_{3}}\Phi_{k}\right)-2\alpha D_{x_{3}}I_{\alpha}(\Phi_{k})+D_{x_{3}}\Psi_{k}^{\ep}\\
	&\quad-\alpha\left[E^{+}\left(\Phi_{k}\right)-2\alpha I_{\alpha}(\Phi_{k})+\Psi_{k}^{\ep}\right]\\
	&=E^{-}\left(D_{x_{3}}\Phi_{k}\right)-\alpha E^{-}(\Phi_{k})+D_{x_{3}}\Psi_{k}^{\ep}-\alpha\Psi_{k}^{\ep}.
\end{align*}
In the second equality, we have used  \eqref{2.5} and $2 \int_{\R^3_{+}}\Phi_{k}(y)\eta_{\ep}(x-y^*)\dy = E^{+}(\Phi_{k})-E^{-}(\Phi_{k})$.
Thus, using $\eta_\ep^-(x',0,y)=0$, we obtain that $\left(D_{x_{3}}\Phi_{k}^{\ep}-\alpha\Phi_{k}^{\ep}\right)|_{x_{3}=0}=0$ for $k=1,2$.
We also have $\Phi_{3}^{\ep}|_{x_{3}=0}=0$ using \eqref{2.13}, the definition of $\Psi_{3}^{\ep}$. Hence, $\Phi^{\ep}$ satisfy the boundary condition \eqref{NavierBC3}.

\medskip

$\bullet$ We claim that
\begin{equation}\label{2.18}
\|\Phi^{\ep}-\Phi\|_{L^{\infty}\left(-1,0;W^{1,q^{\prime}}\left(\C_{+}(1)\right)\right)}+\|\partial_{t}\Phi^{\ep}-\partial_{t}\Phi\|_{L^{q^{\prime},\infty}(Q_{+}(1))}\rightarrow0,
\end{equation}
as $\ep\rightarrow0$.
Actually, by Lemma \ref{lem0a}(ii), we have
\begin{align}\label{2.19}
	\|I_{\alpha}(\Phi)\|_{L^{q^{\prime},\infty}(Q_{+}(1))}+	\|D_{x^{\prime}} I_{\alpha}(\Phi)\|_{L^{q^{\prime},\infty}(Q_{+}(1))}+	\|\partial_{t}I_{\alpha}(\Phi)\|_{L^{q^{\prime},\infty}(Q_{+}(1))}\lesssim \ep^{1+\frac{1}{q^{\prime}}},
\end{align}
and
\begin{align}\label{2.20}
	\|D_{x_{3}} I_{\alpha}(\Phi)\|_{L^{q^{\prime},\infty}(Q_{+}(1))}\lesssim \ep^{\frac{1}{q^{\prime}}}.
\end{align}
Since $\Psi_{3}^{\ep}(x,t)=-2\alpha\phi(x_{3})\int_{0}^{\ep}e^{-\alpha\xi}\int_{\R^3_{+}} \Phi_{3}(y,t)\cdot\eta_{\ep}(x^{\prime}-y^{\prime},y_{3}+\xi)\dy\dxi$, we can obtain that
\begin{equation}\label{2.21}
	\|\Psi_{3}^{\ep}\|_{L^{\infty}\left(-1,0;W^{1,q^{\prime}}\left(\C_{+}(1)\right)\right)}+\|D_{x_{3}}^2 \Psi_{3}^{\ep}\|_{L^{q^{\prime},\infty}(Q_{+}(1))}+	\|\partial_{t}\Psi_{3}^{\ep}\|_{L^{q^{\prime},\infty}(Q_{+}(1))}\lesssim \ep.
\end{equation}
Thus, by \eqref{Bogovskii} and \cite[Lemma \textrm{III}.3.1]{Galdi2011}, we have that for $k=1,2$
\begin{equation}\label{2.22}
\|\Psi_{k}^{\ep}\|_{L^{\infty}\left(-1,0;W^{1,q^{\prime}}\left(\C_{+}(1)\right)\right)}+	\|\partial_{t}\Psi_{k}^{\ep}\|_{L^{q^{\prime},\infty}(Q_{+}(1))}\lesssim \ep.
\end{equation}
Moreover, by \eqref{2.7} and the fact that $\Phi$ and all their derivatives are uniformly continuous in $Q_{+}(1)$, we have
\begin{equation}\label{2.23}
	\|E(\Phi)-\Phi\|_{L^{\infty}\left(-1,0;W^{1,q^{\prime}}\left(\C_{+}(1)\right)\right)}+\|\partial_{t}E(\Phi)-\partial_{t}\Phi\|_{L^{q^{\prime},\infty}(Q_{+}(1))}\rightarrow0,
\end{equation}
as $\ep\rightarrow0$. Combining \eqref{2.19}-\eqref{2.23}, we have \eqref{2.18}.

Since $\bm{u}$ is a very weak solution, testing \eqref{Stokes-weak-form} with the valid test function $\Phi^{\ep}$, we obtain
\begin{align*}
   \iint_{Q_{+}(1)} & \bm{f}\cdot\Phi^{\ep}-\CF:\nabla \Phi ^{\ep} \dx\dt
=	- \iint_{Q_{+}(1)} \bm{u} \cdot \left(\partial_{t}+\Delta\right) \Phi^{\ep} \dx\dt\\
&= \iint_{Q_{+}(1)}  -\bm{u} \cdot \partial_{t}\Phi^{\ep}+\nabla \bm{u}\cdot\nabla  \Phi^{\ep} \dx\dt+\sum_{k=1,2}  \iint_{\Sigma}\alpha u_{k}\cdot\Phi^{\ep}_{k}\dx^{\prime}\dt.
\end{align*}
Sending $\ep\rightarrow0$ and using \eqref{2.18}, we get \eqref{weakform}, and $\bm{u}$ is a weak solution.
\end{proof}

\section{Boundary gradient estimates}\label{Sec 4}
In this section we will prove Theorem \ref{thm1}. The idea is to choose suitable test fields $\Phi$ in the definition of very weak solutions to get uniform gradient bound of approximations of $\bm{u}$. We will define $\bm{w}$ which are approximations of the vorticity $\curl \bm{u}$ plus some lower order terms, and satisfy the nonhomogeneous heat equation in the whole space $\R^3$.

In the appendix \S\ref{sec:app}, we will give an alternative proof of  Theorem \ref{thm1} which uses the boundary condition \eqref{NavierBC3} of $\curl \bm{u}$ to estimate $\curl \bm{u}$ in the half space $\R^3_+$. It is more intuitive, but assumes higher regularity  of $\bm{u}$.

By assumption, $\bm{u},\CF\in L^{q,r}\left(Q_{+}(1)\right)$ and $\bm{f}\in L^{q_{*},r}\left(Q_{+}(1)\right)$ for $1<q,r<\infty$. We first introduce the approximation vector $\bm{w}=\left(w_{1},w_{2},w_{3}\right)\in L^{1}\left(-1,0;C_{c}^{\infty}(\C(1))\right)$  with
\begin{align*}
	w_{1}(x,t)
	&=D_{x_{2}}\int_{\R^3_{+}}e^{\alpha y_{3}/2 }\zeta(y)u_{3}(y)\eta^-_\ep (x,y)\dy-D_{x_{3}}\int_{\R^3_{+}}e^{\alpha y_{3}/2 }\zeta(y) u_{2}(y)\eta^+_\ep(x,y)\dy\\
	&=\int_{\R^3_{+}}e^{\alpha y_{3}/2 }\zeta(y)\left(u_{2}D_{y_{3}}-u_{3}D_{y_{2}}\right)\eta^{-}_{\ep}(x,y)\dy\\
	&=\int_{\R^3_{+}}\left(u_{2}D_{y_{3}}-u_{3}D_{y_{2}}\right)\left(e^{\alpha y_{3}/2}\zeta(y)\eta^{-}_{\ep}(x,y)\right)\dy\\
	&\quad-\int_{\R^3_{+}}\left( u_{2}D_{y_{3}}-u_{3}D_{y_{2}}\right)\left(e^{\alpha y_{3}/2}\zeta(y)\right)\cdot\eta^{-}_\ep(x,y)\dy
	=:w_{1}^{H}-w_{1}^{L}.
\end{align*}
We have used Lemma \ref{lem0}(ii) in the first equality. We
similarly define
\begin{align*}
	w_{2}(x,t)&=D_{x_{3}}\int_{\R^3_{+}}e^{\alpha y_{3}/2 }\zeta(y)u_{1}(y)\eta^+_\ep (x,y)\dy-D_{x_{1}}\int_{\R^3_{+}}e^{\alpha y_{3}/2 }\zeta(y) u_{3}(y)\eta^-_\ep(x,y)\dy\\
	&=\int_{\R^3_{+}}\left(u_{3}D_{y_{1}}-u_{1}D_{y_{3}}\right)\left(e^{\alpha y_{3}/2}\zeta(y)\eta^{-}_{\ep}(x,y)\right)\dy\\
	&\quad-\int_{\R^3_{+}} \left(u_{3}D_{y_{1}}-u_{1}D_{y_{3}}\right)\left(e^{\alpha y_{3}/2}\zeta(y)\right)\cdot \eta^{-}_{\ep}(x,y)\dy
	=:w_{2}^{H}-w_{2}^{L},
\end{align*}
and
\begin{align*}
	w_{3}(x,t)&=D_{x_{1}}\int_{\R^3_{+}}e^{\alpha y_{3} }\zeta(y)u_{2}(y)\eta^+_\ep (x,y)\dy-D_{x_{2}}\int_{\R^3_{+}}e^{\alpha y_{3} }\zeta(y) u_{1}(y)\eta^+_\ep(x,y)\dy\\
	&=\int_{\R^3_{+}}\left(u_{1}D_{y_{2}}-u_{2}D_{y_{1}}\right)\left(e^{\alpha y_{3}}\zeta(y)\eta^{+}_{\ep}(x,y)\right)\dy\\
	&\quad-\int_{\R^3_{+}} \left(u_{1}D_{y_{2}}-u_{2}D_{y_{1}}\right)\left(e^{\alpha y_{3}}\zeta(y)\right)\cdot\eta^{+}_\ep(x,y)\dy
	=:w_{3}^{H}-w_{3}^{L}.
\end{align*}
Note that the exponential factor $e^{\alpha y_{3} }$ in the integrand of $w_3$ differs from those for $w_1$ and $w_2$. The choices of these factors are to ensure the boundary condition \eqref{NavierBC3} of test fields \eqref{Phi1} and those to be used before \eqref{eq3.4}.  %

Consider the vector field
\begin{equation}\label{Phi1}
\Phi(x,y,t)=\psi(t)\left(0,D_{y_{3}},-D_{y_{2}}\right)\left(\xi(y)\eta^{-}_{\ep}(x,y)\right), \quad \xi(y)=e^{\alpha y_{3}/2}\zeta(y),
\end{equation}
where $\psi\in C^1_c(-1,0)$ is an arbitrary cutoff function. For fixed $x$, it is divergence free and
satisfies the boundary condition \eqref{NavierBC3} at $y_3=0$, using Lemma \ref{lem0}(i), $\eta_\ep^-(x,y)|_{y_3=0}=D_{y_3}^2\eta_\ep^-(x,y)|_{y_3=0}=0$, and the exponent $\al/2$ in $e^{\alpha y_{3}/2}$. Using the valid test field $\Phi$ in the definition of very  weak solution \eqref{Stokes-weak-form}, we can show
\begin{align}
	&\quad\partial_{t}w_{1}^{H}-\Delta w_{1}^{H} \notag\\
	&=\int_{\R^3_{+}}\left(u_{2}D_{y_{3}}-u_{3}D_{y_{2}}\right)\left(2\nabla_{y}\left(\xi(y)\right)\cdot \nabla_{y}\eta_{\ep}^{-}(x,y)+\Delta_{y} \left(\xi(y)\right)\eta_{\ep}^{-}(x,y)\right)\dy \label{0529a}\\
	&\quad+\int_{\R^3_{+}}\left(f_{2}D_{y_{3}}-f_{3}D_{y_{2}}-\tsum_{k=1}^3\left(F_{2k}D_{y_{3}}-F_{3k}D_{y_2}\right)D_{y_{k}}\right)\left(\xi(y)\eta^{-}_{\ep}(x,y)\right)\dy. \notag
\end{align}
Indeed, writing 
\[
w_1^H(x,t) = \int_{\R^3_{+}} L (\xi \eta^-_\ep)\dy,\quad L = u_2D_{y_{3}}-u_{3}D_{y_{2}},
\]
we have
\begin{align*}
	-\Delta_{x} w_{1}^{H}=&-\int_{\R^3_{+}}L\left[\xi(y) \Delta_{x}\eta^{-}_{\ep}(x,y)\right]\dy
	=-\int_{\R^3_{+}}L \left[\xi(y) \Delta_{y}\eta^{-}_{\ep}(x,y)\right]\dy\\
	=&-\int_{\R^3_{+}}L \Delta_{y}(\xi(y) \eta^{-}_{\ep}(x,y))\dy + I_1,
\end{align*}
where
\[	
I_1	=\int_{\R^3_{+}}L \left[2\nabla_{y}\left(\xi(y)\right)\cdot \nabla_{y}\eta_{\ep}^{-}(x,y)+\left(\Delta_{y} \xi(y)\right)\eta_{\ep}^{-}(x,y)\right]\dy .
\]	
Hence for $\psi(t) \in C^1_c(-1,0)$,
\EQN{
&\int_{-1}^0 \psi(t)(\pd_t -\De_x) w_1^H(x,t) \dt
\\
&=\int_{-1}^0-w_1^H \psi'(t)\dt + \int_{-1}^0 \psi(t)\bke{-\int_{\R^3_{+}}L \Delta_{y}(\xi \eta^{-}_{\ep})\dy + I_1}\dt
\\
&=-\int_{-1}^0 \int_{\R^3_{+}} L (\pd_t+\Delta_{y})[\psi(t)\xi \eta^{-}_{\ep}]\dy\dt + \int_{-1}^0 \psi(t)I_1\dt.
}

From the definition \eqref{Stokes-weak-form} with $\Phi = \psi(t)(0,D_{y_3}, -D_{y_2})(\xi(y) \eta^{-}_{\ep})$, we get
\[
-\int_{-1}^0 \int_{\R^3_{+}}  L (\pd_t+\Delta_{y})[\psi(t)\xi \eta^{-}_{\ep}]\dy\dt=\int_{-1}^0\psi (t)I_F \dt,
\]
where
\[
I_F=\int_{\R^3_{+}}\left[ f_{2}D_{y_{3}}-f_{3}D_{y_{2}}-\tsum_{k=1}^3\left(F_{2k}D_{y_{3}}-F_{3k}D_{y_2}\right)D_{y_{k}}\right]\left(\xi(y)\eta^{-}_{\ep}(x,y)\right)\dy.
\]
Hence we conclude
for any $\psi \in C^1_c(-1,0)$,
\EQN{
&\int_{-1}^0 \psi(t) (\pd_t -\De_x) w_1^H(x,t) \dt
=\int_{-1}^0 \psi(t)I_F\dt + \int_{-1}^0 \psi(t)I_1\dt,
}
or
$
(\pd_t -\De_x) w_1^H(x,t) = I_F + I_1$, which is exactly \eqref{0529a}.

To estimate $w_1^H(x,t)$ using \eqref{0529a}, we change $D_y$ derivatives acting on $\eta_\ep^-(x,y)$ to $D_x$ derivatives using Lemma \ref{lem0} (ii),
and pull them outside of the integrals. For example, 
\EQS{\label{eq3.3}
&\int_{\R^3_{+}}u_{2}D_{y_{3}}\left(\nabla_{y}\left(\xi(y)\right)\cdot \nabla_{y}\eta_{\varepsilon}^{-}(x,y)\right)\dy\\
=&\tsum_k\int_{\R^3_{+}}u_{2} \bkt{ D_{y_{3}} D_{y_k}\xi(y) D_{y_k} \eta_{\varepsilon}^{-}
+	D_{y_k}\xi(y) D_{y_{3}} D_{y_k} \eta_{\varepsilon}^{-}} \dy\\
=&-\tsum_k D_{x_{k}}\int_{\R^3_{+}}u_{2}\left(D_{{3}}D_{{k}}\xi(y)\right)\cdot\eta_{\ep}^{-\epsilon _k}(x,y)\dy\\
&+\tsum_k D_{x_{k}}D_{x_{3}}\int_{\R^3_{+}}u_{2}\left(D_{{k}}\xi(y)\right)\cdot\eta_{\ep}^{+\epsilon _k}(x,y)\dy.
}
Other terms are similar. Furthermore, we may cut off the initial data of $w_1^H$ by considering $(\pd_t-\Delta) (\th w_1^H)=\th (\pd_t-\Delta) (w_1^H)+ \th' w_1^H$ with the error term $\th' w_1^H$ written as similar integrals.
By standard estimates of the heat equation in $\R^3$, Sobolev embedding inequality, Lemma \ref{lem0} (iv), and the fact that $w_1^H$ has compact spatial support, we have
\begin{align}\label{wH1}
	\|w_{1}^{H}\cdot \theta(t)\|_{L^{q,r}(Q_{+}(1))}\lesssim \|\bm{u}\|_{L^{q,r}(Q_{+}(1))}+\|\bm{f}\|_{L^{q_{*},r}(Q_{+}(1))}+\|\CF\|_{L^{q,r}(Q_{+}(1))}.
\end{align}
(For the part involving $\bm{f}$, we use potential estimate when $q \le 3/2$ and need $q_*>1$ when $q=3/2$. We use singular integral estimate when $q>3/2$ and $q_*>1$.)

Analogously, with the divergence-free test field $\Phi=\psi(t)\left(-D_{y_{3}},0,D_{y_{1}}\right)\left(\xi(y)\eta^{-}_{\ep}(x,y)\right)$ satisfying the boundary condition \eqref{NavierBC3}, we have
\begin{align*}
	&\quad\partial_{t}w_{2}^{H}-\Delta w_{2}^{H}\\
	&=\int_{\R^3_{+}}\left(u_{3}D_{y_{1}}-u_{1}D_{y_{3}}\right)\left(2\nabla_{y}\left(\xi(y)\right)\cdot \nabla_{y}\eta_{\ep}^{-}(x,y)+\Delta_{y} \left(\xi(y)\right)\eta_{\ep}^{-}(x,y)\right)\dy\\
	&\quad+\int_{\R^3_{+}}\left(f_{3}D_{y_{1}}-f_{1}D_{y_{3}}-\tsum_{k=1}^3\left(-F_{1k}D_{y_{3}}+F_{3k}D_{y_1}\right)D_{y_{k}}\right)\left(\xi(y)\eta^{-}_{\ep}(x,y)\right)\dy,
\end{align*}
and for the test field $\Phi=\psi(t)\left(D_{y_{2}},-D_{y_{1}},0\right)\left(e^{\alpha y_{3}}\zeta(y)\eta^{+}_{\ep}(x,y)\right)$,
\begin{align*}
	&\quad\partial_{t}w_{3}^{H}-\Delta w_{3}^{H}\\
	&=\int_{\R^3_{+}}\left(u_{1}D_{y_{2}}-u_{2}D_{y_{1}}\right)\left(2\nabla_{y}\left(e^{\alpha y_{3}}\zeta(y)\right)\cdot \nabla_{y}\eta_{\ep}^{-}(x,y)+\Delta_{y} \left(e^{\alpha y_{3}}\zeta(y)\right)\eta_{\ep}^{-}(x,y)\right)\dy\\
	&\quad+\int_{\R^3_{+}}\left(f_{1}D_{y_{2}}-f_{2}D_{y_{1}}-\tsum_{k=1}^3\left(F_{1k}D_{y_{2}}-F_{2k}D_{y_1}\right)D_{y_{k}}\right)\left(e^{\alpha y_{3}}\zeta(y)\eta^{-}_{\ep}(x,y)\right)\dy.
\end{align*}
Similar to \eqref{wH1}, we have
\begin{equation}\label{eq3.4}
	\|\bm{w}^{H}\cdot \theta(t)\|_{L^{q,r}(Q_{+}(1))}\lesssim \|\bm{u}\|_{L^{q,r}(Q_{+}(1))}+\|\bm{f}\|_{L^{q_{*},r}(Q_{+}(1))}+\|\CF\|_{L^{q,r}(Q_{+}(1))},
\end{equation}
while
\begin{equation}
	\|\bm{w}^{L}\cdot \theta(t)\|_{L^{q,r}(Q_{+}(1))}\lesssim\|\bm{u}\|_{L^{q,r}(Q_{+}(1))}.
\end{equation}
Accordingly, we get
\begin{align}\label{eq3.6}
	\|\bm{w}\cdot \theta(t)\|_{L^{q,r}(Q_{+}(1))}\lesssim \|\bm{u}\|_{L^{q,r}(Q_{+}(1))}+\|\bm{f}\|_{L^{q_{*},r}(Q_{+}(1))}+\|\CF\|_{L^{q,r}(Q_{+}(1))}.
\end{align}

Denote the vector $\bm{v}=\left(E^{+}\left(u_{1}\zeta\right),E^{+}\left(u_{2}\zeta\right),E^{-}(u_{3}\zeta)\right)$. Now we consider the error $\left(e^{-\alpha x_{3}/2}w_{1},e^{-\alpha x_{3}/2}w_{2},e^{-\alpha x_{3}}w_{3}\right)-\nabla\times \bm{v}$ for $x\in \C_+(1)$. For instance, its first component
\begin{align*}
&\text{Err}_1=e^{-\frac\alpha2 x_{3}}w_{1}-\left(D_{x_{2}}v_{3}-D_{x_{3}}v_{2}\right)=D_{x_{2}}\int_{\R^3_{+}}\left(e^{-\frac\alpha2 (x_{3}-y_{3}) }-1\right)\zeta(y)u_{3}(y)\eta^-_\ep (x,y)\dy\\
	&-D_{x_{3}}\int_{\R^3_{+}}\left(e^{-\frac\alpha2 (x_{3}-y_{3})}-1\right)\zeta(y) u_{2}(y)\eta^+_\ep(x,y)\dy
	-\frac{\alpha}{2}\int_{\R^3_{+}}e^{-\frac\alpha2 (x_{3}-y_{3}) }\zeta(y) u_{2}(y)\eta^+_\ep(x,y)\dy .
\end{align*}
Hence
\begin{align*}
	|\text{Err}_1|
	\lesssim&\int_{\R^3_{+}}|x_{3}-y_{3}|\zeta(y)|\bm{u}(y)||\nabla_{y} \eta_{\ep}^{-}(x,y)|+\zeta(y)|\bm{u}(y)||\eta_{\ep}^{+}(x,y)|\dy \\
	\lesssim& \int_{\R^3_{+}}\zeta(y)|\bm{u}(y)|\left(\varepsilon|\nabla_{y} \eta_{\ep}^{-}(x,y)|+|\eta_{\ep}^{+}(x,y)|\right)\dy.
\end{align*}
Accordingly,
\begin{align*}
\|e^{-\alpha x_{3}/2}w_{1}-\left(D_{x_{2}}v_{3}-D_{x_{3}}v_{2}\right)\|_{L^{q}\left(\C_+(1)\right)}\lesssim \|\bm{u}\|_{L^{q}\left(\C_+(1)\right)}.
\end{align*}
Similarly, we have control of the whole error term
\begin{align*}
	\|\left(e^{-\alpha x_{3}/2}w_{1},e^{-\alpha x_{3}/2}w_{2},e^{-\alpha x_{3}}w_{3}\right)-\nabla\times \bm{v}\|_{L^{q}\left(\C_+(1)\right)}\lesssim \|\bm{u}\|_{L^{q}\left(\C_+(1)\right)}.
\end{align*}
Hence we get from \eqref{eq3.6}
\begin{equation}\label{eq3.8}
	\|\nabla\times \bm{v}\cdot \theta(t)\|_{L^{q,r}(Q_{+}(1))}\lesssim \|\bm{u}\|_{L^{q,r}(Q_{+}(1))}+\|\bm{f}\|_{L^{q_{*},r}(Q_{+}(1))}+\|\CF\|_{L^{q,r}(Q_{+}(1))}.
\end{equation}

Notice that
$\div \bm{v} = \int_{\R^3_+} \bm{u} \cdot \nabla \zeta \eta_\ep^-(x,y)\dy$ and
$\bm{v}\cdot\bm{n}=0$ on the boundary of $\C_+(1)$.
By Theorem 3.5 in \cite[p.55]{Amrouche2012}, we have the uniform-in-$\ep$ estimate
\begin{equation}\label{1order1}
	\|\nabla \bm{v} \cdot\theta(t)\|_{L^{q,r}(Q_{+}(1))}\lesssim\|\bm{u}\|_{L^{q,r}(Q_{+}(1))}+\|\bm{f}\|_{L^{q_{*},r}(Q_{+}(1))}+\|\CF\|_{L^{q,r}(Q_{+}(1))}.
\end{equation}
By Lemma \ref{lem0} (iii), we let $\ep\rightarrow0$ and get $\nabla (\bm{u}\zeta )\cdot\theta(t)\in L^{q,r}(Q_{+}(1))$ with
\begin{equation}\label{0530}
\|\nabla \left(\bm{u}\zeta\right) \cdot\theta(t)\|_{L^{q,r}(Q_{+}(1))}\lesssim \|\bm{u}\|_{L^{q,r}(Q_{+}(1))}+\|\bm{f}\|_{L^{q_{*},r}(Q_{+}(1))}+\|\CF\|_{L^{q,r}(Q_{+}(1))}.
\end{equation}
Hence \eqref{main1}. 

By Lemma \ref{lem2}, a very weak solution $\bm{u}$ satisfying \eqref{main1} is a weak solution. This completes the proof of Theorem \ref{thm1}.

\section{Higher derivative estimates}\label{Sec 5}
In this section we will prove Theorem \ref{thm2}, on the estimates of higher derivatives.
Notice that a weak solution $\bm{u}$ is also a very weak solution and satisfies \eqref{0530}. Moreover, replacing the spatial cutoff function $\zeta$ by its derivative $D_{x}^{\sigma}\zeta$ in the proof of Theorem \ref{thm1}, we obtain that for all $|\sigma|\geq 0$,
\begin{align}\label{1order2}
	\|\nabla \bm{u}\cdot D_{x}^{\sigma}\zeta \cdot\theta(t)\|_{L^{q,r}(Q_{+}(1))}\lesssim \|\bm{u}\|_{L^{q,r}(Q_{+}(1))}+\|\bm{f}+\div\CF\|_{L^{q,r}(Q_{+}(1))}.
\end{align}
Above we have used $q_* \le q$. We also integrated by parts terms involving $\CF$ in the equation of $\bm{w}^{H}$, for instance \eqref{0529a}, by the assumption  $\CF_{13}=\CF_{23}=0$ on the boundary $\Sigma$ and $\eta^-_\ep(x,y',0)=0$. 
Estimate \eqref{1order2} will be used in subsections 4.1 and 4.2.
This new cut off function $D_{x}^{\sigma}\zeta$ will be used in subsection 4.3 to bound the higher derivatives.

\subsection{Gradient estimate of the pressure} 
We claim that for all $|\sigma|\geq 0$,
\begin{align}\label{1order-pressure}
	\|\nabla \pi\cdot D_{x}^{\sigma}\zeta \cdot\theta(t)\|_{L^{q,r}(Q_{+}(1))}\lesssim \|\bm{u}\|_{L^{q,r}(Q_{+}(1))}+\|\pi\|_{L^{q,r}(Q_{+}(1))}+\|\bm{f}+\div\CF\|_{L^{q,r}(Q_{+}(1))}.
\end{align}
Without loss of generality, we may assume that $\sigma=0$.
Fix $x\in Q(1)$ with $x_3 \in (-1,1)$. Taking the test function $ \Phi=\psi(t)\nabla _y\left(\zeta(y)\eta^+_{\ep}(x,y)\right) \in C^1_c(Q_+(1)\cup \Sigma)$ in \eqref{weakform2} for arbitrary $\psi \in C^1_c(-1,0)$, and noting that $\Phi_3|_{y_3=0}=D_{y_3}(\zeta \eta^+_{\ep})|_{y_3=0}=0$, we get
\begin{align*}
-\Delta_{x}&\int _{-1}^0\int_{\R^3_+}\pi(y)\zeta(y)\eta^+_{\ep}(x,y) \psi(t)\dy\dt
\\
&=\int _{-1}^0\int_{\R^3_+}\pi(y)\left[2\nabla_{y}\zeta(y)\cdot\nabla_{y}\eta^+_{\ep}(x,y)+(\Delta_{y}\zeta(y))\eta^+_{\ep}(x,y)\right] \psi(t)\dy\dt\\
&\quad-\alpha\tsum_{k=1}^2\int _{-1}^0\int_{\R^2} u_{k}(y^{\prime},0)D_{y_{k}}\left(\zeta(y^{\prime},0)\eta^+_{\ep}(x,y^{\prime},0)\right)\psi(t)\dy^{\prime} \dt\\
&\quad+\int _{-1}^0\int_{\R^3_+} \bkt{\bm{f}\cdot\nabla_{y}\left(\zeta(y)\eta^+_{\ep}(x,y)\right)-\CF:\nabla^2\left(\zeta(y)\eta^+_{\ep}(x,y)\right)}\psi(t)\dy \dt.
\end{align*}
In the calculation we have used $\int _{-1}^0\int_{\R^3_+}
\bket{ \bm{u}\cdot \Phi \psi'(t) - \nabla \bm{u}: \nabla \Phi \psi(t)}\dy \dt=0$ by integration by parts using 
$\div \bm{u}=0$ and $u_3|_{y_3=0}=\Phi_3|_{y_3=0}=0$. Since $\psi(t) \in C^1_c(-1,0)$ is arbitrary, we can remove the time integration and get
\begin{align}\label{4.3}
-\Delta_{x}\int_{\R^3_+}\pi(y)\zeta(y)\eta^+_{\ep}(x,y)\dy&= I_\pi+I_u + I_f,
\end{align}
where
\begin{align*}
I_\pi &= \int_{\R^3_+}\pi(y)\left[ 2\nabla_{y}\zeta(y)\cdot\nabla_{y}\eta^+_{\ep}(x,y)+(\Delta_{y}\zeta(y))\eta^+_{\ep}(x,y)\right]\dy\\
I_f &=\int_{\R^3_+}\left(\bm{f}+\div\CF\right)\cdot\nabla_{y}\left(\zeta(y)\eta^+_{\ep}(x,y)\right)\dy
\end{align*}
(using $F_{13}=F_{23}=\Phi_3=0$ on the boundary $\Sigma$), and
\begin{align*}
I_u 
=& -\alpha\tsum_{k=1}^2\int_{\R^2} u_{k}(y^{\prime},0)D_{y_{k}}\left(\zeta(y^{\prime},0)\eta^+_{\ep}(x,y^{\prime},0)\right)\dy^{\prime}\\
=&\alpha\tsum_{k=1}^2\int_{\R^3_+}D_{y_{3}} \bkt{u_{k}(y)D_{y_{k}}\left(\zeta(y)\eta^+_{\ep}(x,y)\right)}\dy 
\\
=&\alpha\tsum_{k=1}^2\int_{\R^3_+}\bkt{\bke{D_{y_{3}} u_{k}(y)} D_{y_{k}}  - \bke{D_{y_{k}} u_{k}(y)} D_{y_{3}}}
\left(\zeta(y)\eta^+_{\ep}(x,y)\right)\dy .
\end{align*}

Direct estimate would give a factor $1/\ep$ from $D_y \eta_\ep^\pm$. We use Lemma \ref{lem0} (ii) to change $D_y$ acting on $ \eta_\ep^\pm$ to $D_x$ and pull it outside the integrals, as what we did in \eqref{eq3.3}. Then \eqref{4.3} is of the form $-\De \Pi = f+ \div g$. 
By standard estimates of Laplace equation \eqref{4.3} in $\R^3$ (not $\R^3_+$), Sobolev embedding inequality,  and \eqref{1order2}, we have
\begin{align*}
&\|\nabla_x \int_{\R^3_+}\pi(y)\zeta(y)\eta^+_{\ep}(x,y)\dy\cdot \theta(t)\|_{L^{q,r}(Q(1))}
\\
&\lesssim \|\bm{u}\|_{L^{q,r}(Q_{+}(1))}+\|\pi\|_{L^{q,r}(Q_{+}(1))}
+\|\bm{f}+\div\CF\|_{L^{q,r}(Q_{+}(1))}.
\end{align*} 
Letting $\varepsilon \rightarrow 0_+$ and using Lemma \ref{lem0} (ii)-(iv), we get
\begin{equation*}
\|\nabla (\pi \zeta) \theta(t)\|_{L^{q,r}(Q_{+}(1))}\lesssim \|\bm{u}\|_{L^{q,r}(Q_{+}(1))}+\|\pi\|_{L^{q,r}(Q_{+}(1))}
+\|\bm{f}+\div\CF\|_{L^{q,r}(Q_{+}(1))}.
\end{equation*}

\subsection{Estimate of second derivative of the velocity}\label{sec4.2}
We claim that
\begin{align}\label{2order}
&\quad\|\nabla^2 \bm{u}\cdot D_{x}^{\sigma} \zeta\cdot\theta(t)\|_{L^{q,r}(Q_{+}(1))} +\|\partial_{t}\bm{u}\cdot D_{x}^{\sigma}\zeta\cdot\theta(t)\|_{L^{q,r}(Q_{+}(1))}\notag\\
&\lesssim \|\bm{u}\|_{L^{q,r}(Q_{+}(1))}+\|\pi\|_{L^{q,r}(Q_{+}(1))}
+\|\bm{f}+\div\CF\|_{L^{q,r}(Q_{+}(1))}.
\end{align}
Without loss of generality, we may assume that $\sigma=0$. Denote the vector $\bm{v}=\left(v_{1},v_{2},v_{3}\right)$ with
\begin{align*}
v_{i}(x,t)&=\int_{\R^3_{+}}u_{i}(y,t)e^{-\alpha y_3}\zeta(y)\cdot \eta^{+}_{\ep}(x,y)\dy,\quad (i=1,2),\\
v_{3}(x,t)&=\int_{\R^3_{+}}u_{3}(y,t)\zeta(y)\cdot \eta^{-}_{\ep}(x,y)\dy.
\end{align*}
Using Lemma \ref{lem0}(ii),
integration by parts, and $\pd_{y_3} \eta_{\ep}^{+}(x,y)|_{y_3=0}=0$, we have
\begin{align*}
	-\Delta_{x} v_{1}
	=&-\int_{\R^3_{+}}u_{1}(y,t)e^{-\alpha y_3}\zeta(y)\cdot \Delta_{y}\eta^{+}_{\ep}(x,y)\dy \\
	=&\int_{\R^3_{+}}u_{1}(y,t)\bket{\nabla \cdot  \bkt{ - \nabla_{y}\left(e^{-\alpha y_3}\zeta\eta^{+}_{\ep}\right)
	+ 2 \nabla \left(e^{-\alpha y_3}\zeta\right)\eta^{+}_{\ep}} 
	- \Delta_{y}\left(e^{-\alpha y_{3}}\zeta\right)\eta_{\ep}^{+}	}\dy\\
=&\int_{\R^3_{+}}\nabla_{y} u_{1}\cdot  \nabla_{y}\left(e^{-\alpha y_{3}}\zeta(y)\eta_{\ep}^{+}(x,y)\right)\dy
+\alpha \int_{\R^2}u_{1}(y^{\prime},0,t)\zeta(y^{\prime},0)\eta_{\ep}^{+}(x,y^{\prime},0)\dy^{\prime} 
\\
&\quad-\int_{\R^3_{+}}\left[2 \nabla_{y}u_{1}\cdot\nabla_{y}\left(e^{-\alpha y_{3}}\zeta(y)\right)+u_{1}\Delta_{y}\left(e^{-\alpha y_{3}}\zeta(y)\right)\right]\cdot\eta_{\ep}^{+}(x,y)\dy.
\end{align*}
Taking the test function $\Phi=\psi(t)\left(e^{-\alpha y_{3}}\zeta(y) \eta^+_{\ep}(x,y),0,0\right)$ in \eqref{weakform2} (and removing the cuf-off $\psi(t)$ in $t$ as for  \eqref{0529a}),
we have
\begin{align}\label{v1}
\partial_{t}v_{1}-\Delta_{x} v_{1}&=\int_{\R^3_{+}}\left(f_{1}+\tsum_{k=1}^3 D_{y_{k}}F_{1k}-D_{y_{1}}\pi\right)\cdot e^{-\alpha y_{3}}\zeta(y) \eta^+_{\ep}(x,y)\dy\notag\\
&\quad-\int_{\R^3_{+}}\left[2 \nabla_{y}u_{1}\cdot\nabla_{y}\left(e^{-\alpha y_{3}}\zeta(y)\right)+u_{1}\Delta_{y}\left(e^{-\alpha y_{3}}\zeta(y)\right)\right] \cdot\eta_{\ep}^{+}(x,y)\dy.
\end{align}
Analogously,
\begin{align}\label{v2}
	\partial_{t}v_{2}-\Delta_{x} v_{2}&=\int_{\R^3_{+}}\left(f_{2}+\tsum_{k=1}^3 D_{y_{k}}F_{2k}-D_{y_{2}}\pi\right)\cdot e^{-\alpha y_{3}}\zeta(y) \eta^+_{\ep}(x,y)\dy\notag\\
	&\quad-\int_{\R^3_{+}}\left[2 \nabla_{y}u_{2}\cdot\nabla_{y}\left(e^{-\alpha y_{3}}\zeta(y)\right)+u_{2}\Delta_{y}\left(e^{-\alpha y_{3}}\zeta(y)\right)\right]\cdot\eta_{\ep}^{+}(x,y)\dy,
\end{align}
\begin{align}\label{v3}
	\partial_{t}v_{3}-\Delta_{x} v_{3}&=\int_{\R^3_{+}}\left(f_{3}+\tsum_{k=1}^3 D_{y_{k}}F_{3k}-D_{y_{3}}\pi\right)\cdot \zeta(y) \eta^-_{\ep}(x,y)\dy\notag\\
	&\quad-\int_{\R^3_{+}}\left[2 \nabla_{y}u_{3}\cdot\nabla_{y}\zeta(y)+u_{3}\Delta_{y}\zeta(y)\right]\cdot\eta_{\ep}^{-}(x,y)\dy.
\end{align}
Accordingly, by standard estimates of heat equation in $\R^3$, \eqref{1order2} and \eqref{1order-pressure}, we have 
 \begin{align*}
 	&\quad\|\nabla^2 \bm{v}\cdot\theta(t)\|_{L^{q,r}(Q(1))} +\|\partial_{t}\bm{v}\cdot\theta(t)\|_{L^{q,r}(Q(1))}\\
 	&\lesssim \|\bm{u}\|_{L^{q,r}(Q_{+}(1))}+\|\pi\|_{L^{q,r}(Q_{+}(1))}
 	+\|\bm{f}+\div\CF\|_{L^{q,r}(Q_{+}(1))}.
 \end{align*}
Sending $\ep\rightarrow0_+$, 
we get estimates of $D_x^2$ and $\pd_t$ of $u_i e^{-\al y_3}\zeta$, $i=1,2$, and $u_3\zeta$, and hence
\eqref{2order}.

Combining \eqref{1order-pressure} and \eqref{2order}, we have shown \eqref{main2}.

\subsection{Estimate of third derivative of the velocity}
By standard estimates of the equations \eqref{4.3}, \eqref{v1}-\eqref{v3}, and the estimate  \eqref{1order2}, \eqref{1order-pressure} and \eqref{2order}, we have for $k=1,2$
 \EQS{\label{4.9}
	&\quad\|\nabla^3 u_{k}\zeta\,\theta(t)\|_{L^{q,r}(Q_{+}(1))} +\|\nabla\partial_{t}u_{k}\zeta\,\theta(t)\|_{L^{q,r}(Q_{+}(1))}+\|D_{x_{k}}\nabla \pi \zeta\,\theta(t)\|_{L^{q,r}(Q_{+}(1))}\\
	&\quad+\|D_{x_{k}}\nabla^2 u_{3}\zeta\,\theta(t)\|_{L^{q,r}(Q_{+}(1))} +\|D_{x_{k}}\partial_{t}u_{3}\zeta\,\theta(t)\|_{L^{q,r}(Q_{+}(1))}\\
	&\lesssim \|\bm{u}\|_{L^{q,r}(Q_{+}(1))}+\|\pi\|_{L^{q,r}(Q_{+}(1))}+\|\bm{f}+\div\CF\|_{L^{r}\left(-1,0;W^{1,q}(\C_{+}(1))\right)}.
}
The difficult terms are
\EQ{\label{eq4.10}
D_{x_{3}}^3 u_{3}, \quad D_{x_{3}}\partial_{t}u_{3},\quad \text{and}\quad D_{x_{3}}^2 \pi.
}
 For instance, in the $v_{3}$-equation \eqref{v3}, the source terms contain $\eta_{\ep}^{-}$. Taking derivative $D_{x_{3}}$ on both side of the equation \eqref{v3}, we can not deduce what we want after integrating by parts of the source term $D_{y_{3}} \pi\cdot\zeta(y)D_{y_{3}}\eta_{\ep}^{+}$ since we have no boundary condition of $\pi$.
The estimate of terms in \eqref{eq4.10} %
rely on the divergence-free property of the velocity and the equation of $u_{3}$ in \eqref{NS-Eqn}. Noticing that
\begin{align*}
D_{x_{3}}^3 u_{3}=-D_{x_{3}}^2 \left(D_{x_{1}}u_{1}+ D_{x_2}u_{2}\right),\qquad
D_{x_{3}}\partial_{t}u_{3}=-\partial_{t}\left(D_{x_{1}}u_{1}+ D_{x_2}u_{2}\right),
\end{align*}
and $D_{x_{3}}^2 \pi=D_{x_{3}}\left(f_{3}+\tsum_{k=1}^3D_{x_{k}}F_{3k}-\partial_{t}u_{3}+\Delta u_{3}\right)$, we have
 \EQS{\label{4.10}
	&\|D_{x_{3}}^3 u_{3}\zeta\,\theta(t)\|_{L^{q,r}(Q_{+}(1))} +\|D_{x_{3}}\partial_{t}u_{3}\zeta\,\theta(t)\|_{L^{q,r}(Q_{+}(1))}+\|D_{x_{3}}^2 \pi\zeta\,\theta(t)\|_{L^{q,r}(Q_{+}(1))}\\
	&\lesssim \|\bm{u}\|_{L^{q,r}(Q_{+}(1))}+\|\pi\|_{L^{q,r}(Q_{+}(1))}+\|\bm{f}+\div\CF\|_{L^{r}\left(-1,0;W^{1,q}(\C_{+}(1))\right)}.
}
Combing \eqref{4.9} and \eqref{4.10}, we have  \eqref{main3}.

\section{%
Estimates for the case $\alpha=0$}\label{Sec 3}
In this section we prove Theorem \ref{thm0} with $\al=0$, i.e., the infinite slip length case.

We first show a lemma.
\begin{lem} \label{lem3}
	Let the friction coefficient $\alpha=0$. $\bm{u}\in L^{1}(-1,0;W^{1,1}(\C_{+}(1)))$ is a weak solution of the Stokes equations \eqref{Stokes-Eqn} in $Q_+(1)$ with Navier boundary condition  \eqref{NavierBC} on the flat boundary $\Sigma$, and $\left(\overline{\bm{u}},\overline{\bm{f}},\overline{\CF}\right)$ are natural extensions of $\left(\bm{u},\bm{f},\CF\right)(x,t)$ to $Q(1)$: For $(x,t)\in Q_{+}(1)$, let $\left(\overline{\bm{u}},\overline{\bm{f}},\overline{\CF}\right)(x,t)=\left(\bm{u},\bm{f},\CF\right)(x,t)$ and
	\begin{equation}\label{even-odd-ext}
	\bar u_i (x^*,t) = \left \{ \begin{array}{rl} 
	\bar u_i (x,t), \quad  & i=1,2,\\[2pt]
	-\bar u_i (x,t), \quad  & i=3,
\end{array} \right.
	\end{equation}
	\begin{equation}
	\bar f_i (x^*,t) = \left \{ \begin{array}{rl} 
		\bar f_i (x,t), \quad  & i=1,2,\\[2pt]
		-\bar f_i (x,t), \quad  & i=3,
	\end{array} \right.
\end{equation}
	\begin{equation}
	\bar \CF_{ij} (x^*,t) = \left \{ \begin{array}{rl} 
		\bar \CF_{ij} (x,t), \quad  & 1\leq i,j\leq 2\ \text{or}\ i=j=3,\\[2pt]
		-\bar \CF_{ij} (x,t), \quad  & otherwise.
	\end{array} \right.
\end{equation}
Then $\overline{\bm{u}}$ is a weak solution of the Stokes equations \eqref{Stokes-Eqn} in $Q(1)$ with data $\overline{\bm{f}}, \overline{\CF}$, i.e., $\overline{\bm{u}} \in L^{1}\left(-1, 0; W^{1,1}(\C(1))\right)$ is divergence free and satisfies
\begin{equation}\label{weak NS2}
	\iint_{Q(1)} -\overline{\bm{u}} \cdot \partial_{t}\Phi+\nabla \overline{\bm{u}}:\nabla \Phi \dx\dt= \iint_{Q(1)} \overline{\bm{f}}\cdot\Phi-\overline{\CF}:\nabla\Phi\dx\dt
\end{equation}
for all divergence free vector $\Phi \in C^1_{c}(Q(1))$. We do not assert any boundary condition of $\overline{\bm{u}} $ on $\pd Q(1)$.
\end{lem}

\begin{rem}(i)
The even-odd extension \eqref{even-odd-ext} in the 2D setting has also been used in \cite[Section 2]{MR4558989} to extend solutions with Navier BC and $\alpha=0$. As noted in its Remark 2.1, such an extension is not a solution if the Navier BC with $\al=0$ is replaced by the no-slip BC.
\quad
(ii) If we denote $\epsilon_i= 1-2 \de_{i3}$ (i.e., $\epsilon_1=\epsilon_2=1$ and $\epsilon_3=-1$) as in \cite{KLLT1}, then %
\[\bar u_i (x^*,t) = \epsilon_i	\bar u_i (x,t),\quad
\bar f_i (x^*,t) = \epsilon_i	\bar f_i (x,t),\quad
\bar \CF_{ij} (x^*,t) = \epsilon_i \epsilon_j \CF_{ij} (x,t).
\]

\end{rem}

\begin{proof}
It is standard to show that $\div \overline{\bm{u}}=0$ in $Q(1)$ using $u_3=0$ on the boundary $\Sigma$. Since $\bm{u}$ is a weak solution of the Stokes equations \eqref{Stokes-Eqn} in $Q_+(1)$ with $\al=0$, the weak form
\begin{equation}\label{weak NS1}
	\iint_{Q_{+}(1)} -\bm{u} \cdot \partial_{t}\Phi+\nabla \bm{u}:\nabla \Phi \dx\dt= \iint_{Q_{+}(1)} \bm{f}\cdot\Phi-\CF:\nabla\Phi \dx\dt
\end{equation}
holds for all divergence free vector $\Phi \in C^1_{c}(Q_{+}(1)\cup \Sigma)$ with $\Phi_{3}=0$ on the boundary $\Sigma$.
Therefore, by changing variables $x$ to $x^*$, it is easy to prove that \eqref{weak NS2} is valid for all divergence free vector $\Phi \in C^1_{c}(Q(1))$ with $\Phi_{3}=0$ on the boundary $\Sigma$. 

We now start to prove that the weak form \eqref{weak NS2} holds  for all divergence free vector $\Phi \in C^1_{c}(Q(1))$. Decompose $\Phi=\left(\Phi-\Psi\right)+\Psi$, where $\Psi \in C_{c}^{1}(Q(1))$ is a divergence free vector satisfying
\begin{equation}
	\Psi_{3}=\Phi_{3}(x^{\prime},0,t)\phi(|x_{3}|),\quad D_{x_{1}}\Psi_{1}+D_{x_{2}}\Psi_{2}=-D_{x_{3}}\Psi=-\Phi_{3}(x^{\prime},0,t)D_{x_{3}}\phi(|x_{3}|).
\end{equation}
Since
\begin{equation}
	\int_{B^{\prime}(1)}\Phi_{3}(x^{\prime},0,t)\dx^{\prime}=-\int_{\C_{+}(1)}\div \Phi\dx=0,
\end{equation}	
we can define $\Psi_{1},\Psi_{2}$ by the Bogovski\u{\i} formula, similar as \eqref{Bogovskii}. As $\Phi -\Psi$ has zero normal component on the boundary $\Sigma$, $(\Phi_{3} -\Psi_{3})|_{x_3=0}=0$, \eqref{weak NS2} is valid with $\Phi$ replaced by $\Phi -\Psi$. 

We now show \eqref{weak NS2} is valid for $\Phi$ replaced by $\Psi$. In fact $\Psi_{i}(x^*,t)=-\Psi_{i}(x,t)$ for $i=1,2$, $\Psi_{3}(x^*,t)=\Psi_{3}(x,t)$, and $\bar u_{i}\cdot\partial_{t} \Psi_{i}$, $\partial_{j}\bar u_{i}\cdot\partial_{j} \Psi_{i}$, $\bar{f}_{i}\cdot\Psi_{i}$,  $\bar{\CF}_{ij}\cdot\partial_{j}\Psi_{i}$ are all odd functions for $x_{3}\in (-1,1)$. Hence  its integral over $Q(1)$ vanishes. 

As $\Phi=\left(\Phi-\Psi\right)+\Psi$, we have shown \eqref{weak NS2} for any divergence free vector $\Phi \in C^1_{c}(Q(1))$.
\end{proof}

\begin{proof}[Proof of Theorem \ref{thm0}]
By Theorem \ref{thm1}, the solution $\bm{u}$ is a weak solution in $Q_{+}(1/2)$ satisfying that
\begin{equation}
\|\nabla \bm{u}\|_{L^{2}(Q_{+}(1/2))}\lesssim 1+ \|\bm{u}\|_{L^{\infty}(Q_{+}(1))}^2.
\end{equation}
By Lemma \ref{lem3}, its natural extension $\overline{\bm{u}}$ is a weak solution in $Q(1/2)$ and satisfies 
\begin{equation}\label{weak NS3}
	\iint_{Q(1)} -\overline{\bm{u}} \cdot \partial_{t}\Phi+\nabla \overline{\bm{u}}:\nabla \Phi \dx\dt= \iint_{Q(1)} \overline{\bm{u}}\otimes\overline{\bm{u}}:\nabla\Phi\dx\dt
\end{equation}
for all divergence free vector $\Phi \in C^1_{c}(Q(1/2))$. Denote the vector $\bm{v}$  with
\begin{equation*}
	\bm{v}=\int_{\R^3}\overline{\bm{u}}(y,t)\zeta(2y)\cdot \eta_{\ep}(x-y)\dy,
\end{equation*}
and $\bm{w}=\nabla\times\bm{v}$ with
\begin{align}
	\bm{w}
	=&\int_{\R^3}\overline{\bm{u}}\times\nabla_{y}\left(\zeta(2y)\eta_{\ep}(x-y)\right)\dy-\int_{\R^3} \overline{\bm{u}}\times\nabla_{y}\zeta(2y)\cdot\eta_\ep(x-y)\dy\notag\\
	=&\bm{w}^H-\bm{w}^L,
\end{align}
and
\begin{equation}
	\div \bm{v}=\int_{\R^3}\overline{\bm{u}}(y,t)\cdot\nabla_{y}\zeta(2y)\cdot \eta_{\ep}(x-y)\dy.
\end{equation}
By the definition of weak solution \eqref{weak NS3}, we have
\begin{align}\label{3.6}
	\partial_{t}\bm{w}^{H}-\Delta \bm{w}^{H}
	=&\int_{\R^3}\overline{\bm{u}}\times\nabla_{y}\left(2\nabla_{y}\zeta(2y)\cdot \nabla_{y}\eta_{\ep}(x-y)+\Delta_{y} \zeta(2y)\eta_{\ep}(x-y)\right)\dy\notag\\
	&-\int_{\R^3}\left(\overline{\bm{u}}\cdot\nabla_{y}\right)\overline{\bm{u}}\times\nabla_{y}\left(\zeta(2y)\eta_{\ep}(x-y)\right)\dy.
\end{align}
By the mathematical induction, standard estimates of heat equation in $\R^3$, Sobolev embedding inequality, 
the Leibniz rule and Gagliardo-Nirenberg inequalities, we can obtain that for all integer $n\geq0$
\begin{equation}
	\|\nabla^{n}\overline{\bm{u}}\cdot\zeta(2x)\theta(4t)\|_{L^{2}\left(Q(1)\right)}\lesssim 1+ \|\bm{u}\|_{L^{\infty}(Q_{+}(1))}^{n+1}.
\end{equation}
Moreover, it is easy to find that
\begin{equation}
	\|\nabla^{n}\partial_{t} \bm{w}^H\cdot\theta(4t)\|_{L^{2}(\R^3\times(-1,0))}\lesssim 1+\|\bm{u}\|^{n+4}_{L^{\infty}(Q_{+}(1))}.
\end{equation}
Thus, for all $n\geq0$
\begin{equation}
	\| \bm{w}^H\cdot\theta(4t)\|_{L^{\infty}\left(-1,0;W^{n,2}(\R^3)\right)}\lesssim 1+\|\bm{u}\|^{n+4}_{L^{\infty}(Q_{+}(1))}.
\end{equation}
On the other hand, the low order term
\begin{align*}
	&\| \bm{w}^L\cdot\theta(4t)\|_{L^{\infty}\left(-1,0;W^{n,2}(\R^3)\right)}+	\| \div\bm{v}\cdot\theta(4t)\|_{L^{\infty}\left(-1,0;W^{n,2}(\R^3)\right)}\\
	\lesssim&\| \overline{\bm{u}}\times\nabla\zeta(2x)\cdot\theta(4t)\|_{L^{\infty}\left(-1,0;W^{n,2}(\R^3)\right)}+\| \overline{\bm{u}}\cdot\nabla\zeta(2x)\cdot\theta(4t)\|_{L^{\infty}\left(-1,0;W^{n,2}(\R^3)\right)}
\end{align*}
By a simple mathematical induction, we can obtain that for all $n\geq0$
\begin{equation*}
\| \overline{\bm{u}}\cdot\zeta(2x)\theta(4t)\|_{L^{\infty}\left(-1,0;W^{n+1,2}(\R^3)\right)}\lesssim 1+\|\bm{u}\|^{n+4}_{L^{\infty}(Q_{+}(1))}.
\end{equation*}
Hence, $\nabla^n \overline{\bm{u}}\in L^{\infty}(Q(1/4))$ for all $n\geq 1$.
\end{proof}

\section{Appendix: An alternative proof of the gradient estimate}
\label{sec:app}

Our previous proof in Section \ref{Sec 4} shows \eqref{main1} of Theorem \ref{thm1} for any very weak solution. If we assume higher regularity of the solution $\bm{u}$, in particular, the heat equation of $\curl \bm{u}$ and the boundary conditions \eqref{NavierBC} and \eqref{NavierBC3} are meaningful, then 
we have the following alternative and more intuitive proof of \eqref{main1}. For what follows, we assume $\bm{f}=0$ and $\CF=0$ for simplicity.

Recall that the vorticity $\om$ satisfies the heat equation in $Q_+(1)$ with the boundary condition \eqref{NavierBComega}.
Denote the cut-off function $\varphi(t,x)=\theta(t)\zeta(x)$ and   
let $\overline{\bm{u}}=\bm{u}\varphi$ and $\overline{\bm{\omega}}=\bm{\omega} \varphi$. The localized vorticity satisfies the nonhomogeneous heat equation in $\R_{+}^3\times(-1,0)$
\begin{align*}
\partial_{t}\overline{\bm{\omega}}-\Delta \overline{\bm{\omega}}=\bm{g},
\end{align*}
with $\bm{g}=\left(g_{1},g_{2},g_{3}\right)=\bm{\omega}\partial_{t}\varphi-2\nabla\bm{\omega}\cdot\nabla\varphi-\bm{\omega}\Delta\varphi$, and the boundary condition
\begin{align}\label{bC1}
	\overline{\omega}_{1}=-\alpha \overline{u}_{2},\quad\overline{\omega}_{2}=\alpha \overline{u}_{1},\quad\partial_{3}\overline{\omega}_{3}=\alpha\overline{\omega}_{3}.
\end{align}
Denote $K(x,y,t,s)=\Gamma(x-y,t-s)-\Gamma(x^{\prime}-y^{\prime},x_{3}+y_{3},t-s)$, which is the Green function of the heat equation in $\R^3_+$ with zero BC.
 The solution is
\begin{align}
	\overline{\bm{\omega}}=\overline{\bm{\omega}}^{M}+\overline{\bm{\omega}}^{D}
\end{align}
with $\overline{\bm{\omega}}^{M}=\left(\overline{\omega}^{M}_{1},\overline{\omega}^{M}_{2},\overline{\omega}^{M}_{3}\right)$,
\begin{align*}
	\overline{\omega}^{M}_{i}(x,t)=\int_{-1}^{t}\int_{\R^3_{+}}K(x,y,t,s)g_{i}(y,s)\dy\ds,
\end{align*}
for $i=1,2,3$ and $\overline{\bm{\omega}}^{D}=\left(\overline{\omega}^{D}_{1},\overline{\omega}^{D}_{2},\overline{\omega}^{D}_{3}\right)$,
\begin{align*}
	\overline{\omega}^{D}_{1}(x,t)&=2\alpha\int_{-1}^{t}\int_{\R^2}D_{x_{3}}\Gamma(x^{\prime}-y^{\prime},x_{3},t-s)\overline{u}_{2}(y^{\prime},0,s)\dy^{\prime}\ds,\\
	\overline{\omega}^{D}_{2}(x,t)&=-2\alpha\int_{-1}^{t}\int_{\R^2}D_{x_{3}}\Gamma(x^{\prime}-y^{\prime},x_{3},t-s)\overline{u}_{1}(y^{\prime},0,s)\dy^{\prime}\ds,\\
	\overline{\omega}^{D}_{3}(x,t)&=-2\int_{0}^{\infty }e^{-\alpha z}D_{x_{3}}\int_{-1}^{t}\int_{\R^3_{+}}\Gamma(x^{\prime}-y^{\prime},x_{3}+y_{3}+z,t-s) g_{3}(y,s)\dy\ds\dz.
\end{align*}
That $\obo_1$ and $\obo_2$ satisfy the boundary condition  is because $-2 D_{x_{3}}
\Gamma(x^{\prime}-y^{\prime},x_{3},t-s)$ is the Poisson kernel of the heat equation in $\R^3_+$. The formula of $\obo_3$ for the Robin BC $D_3\obo_3=\al \obo_3$ is derived in 
Keller \cite[\S5]{Keller1981}.

\begin{lem}\label{lem1}
For $1<q,r<\infty$, we have
\begin{align}\label{3.2}
	\|\overline{\bm{\omega}}\|_{L^{q,r}\left(Q_{+}(1)\right)}\lesssim \|\bm{u}\|_{L^{q,r}\left(Q_{+}(1)\right)}+\|\overline{\bm{u}}(x^{\prime},0,t)\|_{L^{r}\left(0,1;L^{q}(\R^2)\right)}.
\end{align}
\end{lem}
\begin{proof}

Integrating by parts, we have
\begin{align*}
\overline{\bm{\omega}}^{M}(x,t)&=\int_{-1}^{t}\int_{\R^3_{+}}\left(K(x,y,t,s)\left(\partial_{s}\varphi+\Delta\varphi\right)+2\nabla_{y}K(x,y,t,s)\cdot\nabla_{y}\varphi\right)\cdot\bm{\omega}\dy\ds\\
&=\int_{-1}^{t}\int_{\R^3_{+}}\bm{u}\times\nabla_{y}\left[K(x,y,t,s)\left(\partial_{s}\varphi+\Delta\varphi\right)+2\nabla_{y}K(x,y,t,s)\cdot\nabla_{y}\varphi\right]\dy\ds.
\end{align*}	
By the estimates of parabolic singular integral %
operators, we have that for $1<q,r<\infty$,
\begin{align}\label{3.3}
	\|\overline{\bm{\omega}}^{M}\|_{L^{q,r}\left(Q_{+}(1)\right)}\lesssim \|\bm{u}\|_{L^{q,r}\left(Q_{+}(1)\right)}.
\end{align}
For $i=1,2$, since
\begin{align*}
	|\overline{\omega}_{i}^{D}(x,t)| &\lesssim  \int_{-1}^{t}\frac{1}{t-s}e^{-\frac{x_{3}^2}{8(t-s)}}\int_{\R^2}(4\pi(t-s))^{-1}e^{-\frac{|x^{\prime}-y^{\prime}|^2}{4(t-s)}}\,|\overline{\bm{u}}|(y^{\prime},0,s)\dy^{\prime}\ds
\\
&=  \int_{0}^{\infty}\frac{1}{s}e^{-\frac{x_{3}^2}{8s}}\int_{\R^2}\frac{1}{4\pi s}e^{-\frac{|x^{\prime}-y^{\prime}|^2}{4s}}\,|\overline{\bm{u}}|(y^{\prime},0,t-s)\dy^{\prime}\ds,
\end{align*}
we have
\EQN{
\norm{\overline{\omega}_{i}^{D}(\cdot,t)}_{L^q(\R^3_+)} 
&\lesssim 
\int_{0}^{\infty}\frac{1}{s}\norm{e^{-\frac{x_{3}^2}{8s}}}_{L^q(\R^+)} \norm{\int_{\R^2}\frac{1}{s}e^{-\frac{|x^{\prime}-y^{\prime}|^2}{4s}}\,|\overline{\bm{u}}|(y^{\prime},0,t-s)\dy^{\prime}}_{L^q(\R^2)} \ds \\
&\lec  \int_{0}^{\infty}\frac{1}{s}s^{1/2q} \norm{\overline{\bm{u}}(\cdot,0,t-s)}_{L^q(\R^2)} \ds.
}
Thus, by Hardy-Littlewood inequality,
\begin{align}\label{3.4}
\tsum_{i=1,2}
	\|\overline{\omega}_{i}^{D}\|_{L^{q,r}(Q_{+}(1))}\lesssim \|\overline{\bm{u}}(x^{\prime},0,t)\|_{L^{r}\left(-1,0;L^{q}(\R^2)\right)}.
\end{align}

For $\overline{\omega}^{D}_{3}$, first changing $D_{x_3} \to D_{z_3}$ and then
integrating by parts, we have  
\begin{align*}
	\overline{\omega}^{D}_{3}(x,t)&=2\int_{-1}^{t}\int_{\R^3_{+}}\Gamma(x^{\prime}-y^{\prime},x_{3}+y_{3},t-s) g_{3}(y,s)\dy\ds\\
	&\quad-2\alpha\int_{0}^{\infty }e^{-\alpha z}\int_{-1}^{t}\int_{\R^3_{+}}\Gamma(x^{\prime}-y^{\prime},x_{3}+y_{3}+z,t-s) g_{3}(y,s)\dy\ds\dz
\end{align*}
Analogously to the estimates of $\overline{\omega}_{3}^{M}$, it is easy to deduce that
\begin{align}\label{3.5}
	\|\overline{\omega}_{3}^{D}\|_{L^{q,r}(Q_{+}(1))}\lesssim \|\bm{u}\|_{L^{q,r}(Q_{+}(1))}.
\end{align}	
Summing up \eqref{3.3}, \eqref{3.4} and \eqref{3.5}, we have \eqref{3.2}.
\end{proof}

As $\overline{\bm{u}}=\bm{u}\varphi$ and $\overline{\bm{\omega}}=\bm{\omega} \varphi$, we have $\div \overline{\bm{u}}= \bm{u}\cdot  \nabla \varphi  $, $\curl \overline{\bm{u}} = \overline{\bm{\omega}}+ \nabla \varphi \times\bm{u}$ and $ \overline{\bm{u}}\cdot \bm{n}=0$ on the boundary of $\C_+(1)$.
By Theorem 3.5 in \cite[p.55]{Amrouche2012} and Lemma
 \ref{lem1}, we have %
\begin{align}\label{3.17}
\|D \overline{\bm{u}}\|_{L^{q,r}(Q_{+}(1))}
&\lesssim\|\bm{u}\|_{L^{q,r}(Q_{+}(1))}+ \|\overline{\bm{\omega}}\|_{L^{q,r}\left(Q_{+}(1)\right)}
\\
&\lesssim \|\bm{u}\|_{L^{q,r}\left(Q_{+}(1)\right)}+\|\overline{\bm{u}}(x^{\prime},0,t)\|_{L^{r}\left(0,1;L^{q}(\R^2)\right)}.\nonumber
\end{align}

For fixed $t$,
\[
\int_{\R^2} |\overline{\bm{u}}(x^{\prime},0,t)|^q dx' = \bigg|\int_{\R_+^3} D_3 |\overline{\bm{u}}|^qdx_3\,dx' \bigg|
\le \iint_{\R^3_+} q |D_3\overline{\bm{u}}| \,|\overline{\bm{u}}|^{q-1}\,dx
\le q\norm{D \overline{\bm{u}}}_q \,  \norm{ \overline{\bm{u}}}_q^{q-1}.
\]%
Hence for any $\delta>0$,
\[
\|\overline{\bm{u}}(x^{\prime},0,t)\|_{L^{r}\left(0,1;L^{q}(\R^2)\right)}
\le \delta   \norm{D \overline{\bm{u}}}_{L^{q,r}(Q_{+}(1))}+ C_\delta  \norm{ \overline{\bm{u}}}_{L^{q,r}(Q_{+}(1))}.
\]
By taking $\delta>0$ sufficiently small and
absorbing $\delta   \norm{D \overline{\bm{u}}}_{L^{q,r}(Q_{+}(1))}$ by the left side of \eqref{3.17}, we conclude
\[
\|D \overline{\bm{u}}\|_{L^{q,r}(Q_{+}(1))}
\lesssim \|\bm{u}\|_{L^{q,r}\left(Q_{+}(1)\right)},
\]
which shows \eqref{main1}.

\section*{Acknowledgments}
We warmly thank Professors Reinhard Farwig and \v{S}\'{a}rka Ne\v{c}asov\'{a} for kindly providing us many references, and Prof.~Hongjie Dong for pointing us to the reference \cite{DKP}.
The first author thanks the Department of Mathematics and Pacific Institute for the Mathematical Sciences in University of British Columbia, where part of this work was done.
Hui Chen was supported in part by National Natural Science Foundation of China under the grant [12101556] and China Scholarship Council. 
 The research of both Liang and Tsai was partially supported by Natural Sciences and Engineering Research Council of Canada (NSERC) grants RGPIN-2018-04137 and RGPIN-2023-04534.

\bibliography{NavierBC}
\bibliographystyle{abbrv}

\end{document}